
\magnification = \magstep0
\overfullrule = 2pt
\vsize = 525dd
\hsize = 27cc 
\topskip = 13dd
\hoffset = 1.8cm
\voffset = 1.8cm


\font\bfbig = cmbx10 scaled \magstep2   
\font\partf = cmbx10 scaled \magstep1	

\font\eightrm = cmr10 scaled 800        
\font\sixrm = cmr7 scaled 850
\font\fiverm = cmr5
\font\eighti = cmmi10 scaled 800
\font\sixi = cmmi7 scaled 850
\font\fivei = cmmi5
\font\eightit = cmti10 scaled 800
\font\eightsy = cmsy10 scaled 800
\font\sixsy = cmsy7 scaled 850
\font\fivesy = cmsy5
\font\eightsl = cmsl10 scaled 800
\font\eighttt = cmtt10 scaled 800
\font\eightbf = cmbx10 scaled 800
\font\sixbf = cmbx7 scaled 850
\font\fivebf = cmbx5

\font\fivefk = eufm5                    
\font\sixfk = eufm7 scaled 850
\font\sevenfk = eufm7
\font\eightfk = eufm10 scaled 800
\font\tenfk = eufm10

\newfam\fkfam
	\textfont\fkfam=\tenfk \scriptfont\fkfam=\sevenfk
		\scriptscriptfont\fkfam=\fivefk

\def\eightpoint{%
	\textfont0=\eightrm \scriptfont0=\sixrm \scriptscriptfont0=\fiverm
		\def\rm{\fam0\eightrm}%
	\textfont1=\eighti  \scriptfont1=\sixi  \scriptscriptfont1=\fivei
		\def\oldstyle{\fam1\eighti}%
	\textfont2=\eightsy \scriptfont2=\sixsy \scriptscriptfont2=\fivesy
	\textfont\itfam=\eightit \def\it{\fam\itfam\eightit}%
	\textfont\slfam=\eightsl \def\sl{\fam\slfam\eightsl}%
	\textfont\ttfam=\eighttt \def\tt{\fam\ttfam\eighttt}%
	\textfont\bffam=\eightbf \scriptfont\bffam=\sixbf
		\scriptscriptfont\bffam=\fivebf \def\bf{\fam\bffam\eightbf}%
	\textfont\fkfam=\eightfk \scriptfont\fkfam=\sixfk
		\scriptscriptfont\fkfam=\fivefk
	\rm}
\skewchar\eighti='177\skewchar\sixi='177\skewchar\eightsy='60\skewchar\sixsy='60

\def\fk{\fam\fkfam}
\def\petit{\eightpoint}%


\def\today{\ifcase\month\or
	January\or February\or March\or April\or May\or June\or
	July\or August\or September\or October\or November\or December\fi
	\space\number\day, \number\year}
\def\newline{\hfil\break}
\def\framebox#1{\vbox{\hrule\hbox{\vrule\hskip1pt
	\vbox{\vskip1pt\relax#1\vskip1pt}\hskip1pt\vrule}\hrule}}


\newtoks\RunAuthor\RunAuthor={} \newtoks\RunTitle\RunTitle={}
\def\ShortTitle#1#2{\RunAuthor={#1}\RunTitle={#2}}
\headline={\ifnum\pageno=1{\hfil}\else
	\ifodd\pageno {\petit{\the\RunTitle}\hfil\folio}
	\else {\petit\folio\hfil{\the\RunAuthor}} \fi \fi}
\footline={\hfil}
\long\def\Title#1{\hrule\topglue3truecm
	\noindent{\bfbig#1}\vskip12pt\relax}
\long\def\Author#1#2{\noindent{\bf#1}\vskip6pt%
	\noindent{\petit#2}\vskip10pt%
	\noindent{\petit\today}\vskip32pt\relax}
\long\def\Thanks#1#2{{\parindent=20pt\baselineskip=9pt%
	\footnote{\nobreak${}^{#1}$}{\petit#2\par\vskip-9pt}}}
\long\def\Abstract#1{\noindent{\bf Summary.}\enspace#1\bigskip\relax}
\def\Acknowledgements{\goodbreak\vskip21pt\noindent{\bf Acknowledgements.}%
	\enspace\relax}


\def\PLabel#1{\xdef#1{\nobreak(p.\the\pageno)}}


\newcount\SECNO \SECNO=0
\newcount\SUBSECNO \SUBSECNO=0
\newcount\SUBSUBSECNO \SUBSUBSECNO=0
\def\Part#1{\SECNO=0\SUBSECNO=0\SUBSUBSECNO=0
	\vfill\eject\noindent{\partf Part #1}
	\nobreak\vskip12pt\noindent\kern0pt}
\def\Section#1{\SUBSECNO=0\SUBSUBSECNO=0 \advance\SECNO by 1
	\goodbreak\vskip21pt\leftline{\bf\the\SECNO .\ #1}
	\gdef\Label##1{\xdef##1{\nobreak\the\SECNO}}
	\nobreak\vskip12pt\noindent\kern0pt}
\def\SubSection#1{\SUBSUBSECNO=0 \advance\SUBSECNO by 1
	\goodbreak\vskip21pt\leftline{\it\the\SECNO.\the\SUBSECNO\ #1}
	\gdef\Label##1{\xdef##1{\nobreak\the\SECNO.\the\SUBSECNO}}
	\nobreak\vskip12pt\noindent\kern0pt}
\def\SubSubSection#1{\advance\SUBSUBSECNO by 1
	\goodbreak\vskip21pt\leftline{\rm\the\SECNO.\the\SUBSECNO.\the\SUBSUBSECNO\ #1}
	\gdef\Label##1{\xdef##1{\nobreak\the\SECNO.\the\SUBSECNO.\the\SUBSUBSECNO}}
	\nobreak\vskip12pt\noindent\kern0pt}


\long\def\Definition#1#2{\medbreak\noindent{\bf Definition%
	#1.\enspace}{\it#2}\medbreak\smallskip\relax}
\long\def\Theorem#1#2{\medbreak\noindent{\bf Theorem%
	#1.\enspace}{\it#2}\medbreak\smallskip\relax}
\long\def\Lemma#1#2{\medbreak\noindent{\bf Lemma%
	#1.\enspace}{\it#2}\medbreak\smallskip\relax}

\long\def\Remark#1{\medbreak\noindent{\bf Remark.%
	\enspace}{\it#1}\medbreak\smallskip\relax}


\newcount\FOOTNO \FOOTNO=0
\long\def\Footnote#1{\global\advance\FOOTNO by 1
	{\parindent=20pt\baselineskip=9pt%
	\footnote{\nobreak${}^{\the\FOOTNO)}$}{\petit#1\par\vskip-9pt}%
	}\gdef\Label##1{\xdef##1{\nobreak\the\FOOTNO}}}


\newcount\EQNO \EQNO=0
\def\Eqno{\global\advance\EQNO by 1 \eqno(\the\EQNO)%
	\gdef\Label##1{\xdef##1{\nobreak(\the\EQNO)}}}


\newcount\FIGNO \FIGNO=0
\def\Fcaption#1{\global\advance\FIGNO by 1
	{\petit{\bf Fig. \the\FIGNO.~}#1}
	\gdef\Label##1{\xdef##1{\nobreak\the\FIGNO}}}
\def\Figure#1#2{\medskip\vbox{\centerline{\framebox{#1}}
	\centerline{\Fcaption{#2}}}\medskip\relax}

\newcount\TABNO \TABNO=0
\def\Tcaption#1{\global\advance\TABNO by 1
   {\petit{\bf Table. \the\TABNO.~}#1}
   \gdef\Label##1{\xdef##1{\nobreak\the\TABNO}}}

\newread\epsffilein    
\newif\ifepsffileok    
\newif\ifepsfbbfound   
\newif\ifepsfverbose   
\newif\ifepsfdraft     
\newdimen\epsfxsize    
\newdimen\epsfysize    
\newdimen\epsftsize    
\newdimen\epsfrsize    
\newdimen\epsftmp      
\newdimen\pspoints     
\pspoints=1bp          
\epsfxsize=0pt         
\epsfysize=0pt         
\def\epsfbox#1{\global\def\epsfllx{72}\global\def\epsflly{72}%
   \global\def\epsfurx{540}\global\def\epsfury{720}%
   \def\lbracket{[}\def\testit{#1}\ifx\testit\lbracket
   \let\next=\epsfgetlitbb\else\let\next=\epsfnormal\fi\next{#1}}%
\def\epsfgetlitbb#1#2 #3 #4 #5]#6{\epsfgrab #2 #3 #4 #5 .\\%
   \epsfsetgraph{#6}}%
\def\epsfnormal#1{\epsfgetbb{#1}\epsfsetgraph{#1}}%
\def\epsfgetbb#1{%
%
%
\openin\epsffilein=#1
\ifeof\epsffilein\errmessage{I couldn't open #1, will ignore it}\else
%
%
   {\epsffileoktrue \chardef\other=12
    \def\do##1{\catcode`##1=\other}\dospecials \catcode`\ =10
    \loop
       \read\epsffilein to \epsffileline
       \ifeof\epsffilein\epsffileokfalse\else
%
%
          \expandafter\epsfaux\epsffileline:. \\%
       \fi
   \ifepsffileok\repeat
   \ifepsfbbfound\else
    \ifepsfverbose\message{No bounding box comment in #1; using defaults}\fi\fi
   }\closein\epsffilein\fi}%
%
%
%
\def\epsfclipoff{\def\epsfclipstring{\ifepsfdraft\space clip\fi}}%
\epsfclipoff
\def\epsfsetgraph#1{%
   \epsfrsize=\epsfury\pspoints
   \advance\epsfrsize by-\epsflly\pspoints
   \epsftsize=\epsfurx\pspoints
   \advance\epsftsize by-\epsfllx\pspoints
%
%
   \epsfxsize\epsfsize\epsftsize\epsfrsize
   \ifnum\epsfxsize=0 \ifnum\epsfysize=0
      \epsfxsize=\epsftsize \epsfysize=\epsfrsize
      \epsfrsize=0pt
%
%
     \else\epsftmp=\epsftsize \divide\epsftmp\epsfrsize
       \epsfxsize=\epsfysize \multiply\epsfxsize\epsftmp
       \multiply\epsftmp\epsfrsize \advance\epsftsize-\epsftmp
       \epsftmp=\epsfysize
       \loop \advance\epsftsize\epsftsize \divide\epsftmp 2
       \ifnum\epsftmp>0
          \ifnum\epsftsize<\epsfrsize\else
             \advance\epsftsize-\epsfrsize \advance\epsfxsize\epsftmp \fi
       \repeat
       \epsfrsize=0pt
     \fi
   \else \ifnum\epsfysize=0
     \epsftmp=\epsfrsize \divide\epsftmp\epsftsize
     \epsfysize=\epsfxsize \multiply\epsfysize\epsftmp
     \multiply\epsftmp\epsftsize \advance\epsfrsize-\epsftmp
     \epsftmp=\epsfxsize
     \loop \advance\epsfrsize\epsfrsize \divide\epsftmp 2
     \ifnum\epsftmp>0
        \ifnum\epsfrsize<\epsftsize\else
           \advance\epsfrsize-\epsftsize \advance\epsfysize\epsftmp \fi
     \repeat
     \epsfrsize=0pt
    \else
     \epsfrsize=\epsfysize
    \fi
   \fi
%
%
   \ifepsfverbose\message{#1: width=\the\epsfxsize, height=\the\epsfysize}\fi
   \epsftmp=10\epsfxsize \divide\epsftmp\pspoints
   \vbox to\epsfysize{\vfil\hbox to\epsfxsize{%
      \ifnum\epsfrsize=0\relax
        \includegraphics{\ifepsfdraft}%
      \else
        \epsfrsize=10\epsfysize \divide\epsfrsize\pspoints
        \includegraphics{\ifepsfdraft}%
      \fi
      \hfil}}%
\global\epsfxsize=0pt\global\epsfysize=0pt}%
%
%
{\catcode`\%=12 \global\let\epsfpercent=
%
%
\long\def\epsfaux#1#2:#3\\{\ifx#1\epsfpercent
   \def\testit{#2}\ifx\testit\epsfbblit
      \epsfgrab #3 . . . \\%
      \epsffileokfalse
      \global\epsfbbfoundtrue
   \fi\else\ifx#1\par\else\epsffileokfalse\fi\fi}%
%
%
\def\epsfempty{}%
\def\epsfgrab #1 #2 #3 #4 #5\\{%
\global\def\epsfllx{#1}\ifx\epsfllx\epsfempty
      \epsfgrab #2 #3 #4 #5 .\\\else
   \global\def\epsflly{#2}%
   \global\def\epsfurx{#3}\global\def\epsfury{#4}\fi}%
%
%
\def\epsfsize#1#2{\epsfxsize}
%
%



\newcount\REFNO \REFNO=0
\newbox\REFBOX \setbox\REFBOX=\vbox{}
\def\BegRefs{\setbox\REFBOX\vbox\bgroup
	\parindent18pt\baselineskip9pt\petit}
\def\EndRefs{\par\egroup}
\def\Ref#1{\global\advance\REFNO by 1 \ifnum\REFNO>1\vskip3pt\fi
	\item{\the\REFNO .~}\xdef#1{\nobreak[\the\REFNO]}}
\def\References{\goodbreak\vskip21pt\leftline{\bf References}
	\nobreak\vskip12pt\unvbox\REFBOX\vskip21pt\relax}


\def\N{I\kern-.8ex N}
\def\Z{\raise.72ex\hbox{${}_{\not}$}\kern-.45ex {\rm Z}}
\def\Q{\raise.82ex\hbox{${}_/$}\kern-1.35ex Q} \def\R{I\kern-.8ex R}
\def\C{\raise.87ex\hbox{${}_/$}\kern-1.35ex C} \def\H{I\kern-.8ex H}

\def\RP{{\R\!P}} \def\CP{{\C\!P}} \def\HP{{\H\!P}}

\def\II{{I\!\!I}}

\def\Ct#1{{C\!#1}} \def\Gt#1{{G\!#1}}
\def\Dt#1#2{{D_{\!#1}\!#2}} \def\Tt#1#2{{T_{\!#1}\!#2}}


\Title{M\"obius geometry of surfaces of constant mean curvature
	1 in hyperbolic space}
\Author{U.\ Hertrich-Jeromin\Thanks{\ast}{Partially supported by
	GNSAGA, Italy.},
	E.\ Musso, and L.\ Nicolodi\Thanks{\ast\ast}{Partially supported by
	NSF grant DMS93-12087 and GNSAGA, Italy.}}
	{Dept.\ Mathematics, Technical University Berlin, D-10623 Berlin
	\newline
	Dept.\ Mathematics, Universit\`{a} degli Studi di l'Aquila,
	I-67010 l'Aquila\newline
	Dept.\ Mathematics, Universit\`{a} degli Studi di Roma ``La Sapienza'',
	I-00185 Roma}
\ShortTitle{U.Jeromin, E.Musso, L.Nicolodi}
	{Transformations of Isothermic surfaces}
\Abstract{Various transformations of isothermic surfaces are discussed and
	their interrelations are analyzed. Applications to cmc-1 surfaces
	in hyperbolic space and their minimal cousins in Euclidean space are
	presented: the Umehara-Yamada perturbation, the classical and Bryant's
	Weierstrass type representations, and the duality for cmc-1 surfaces
	are interpreted in terms of transformations of isothermic surfaces.
	A new Weierstrass type representation is introduced and
	a M\"obius geometric characterization of cmc-1 surfaces in hyperbolic
	space and minimal surfaces in Euclidean space is given.}


\def\DPic{\Figure{\epsfxsize.9\hsize\epsfbox{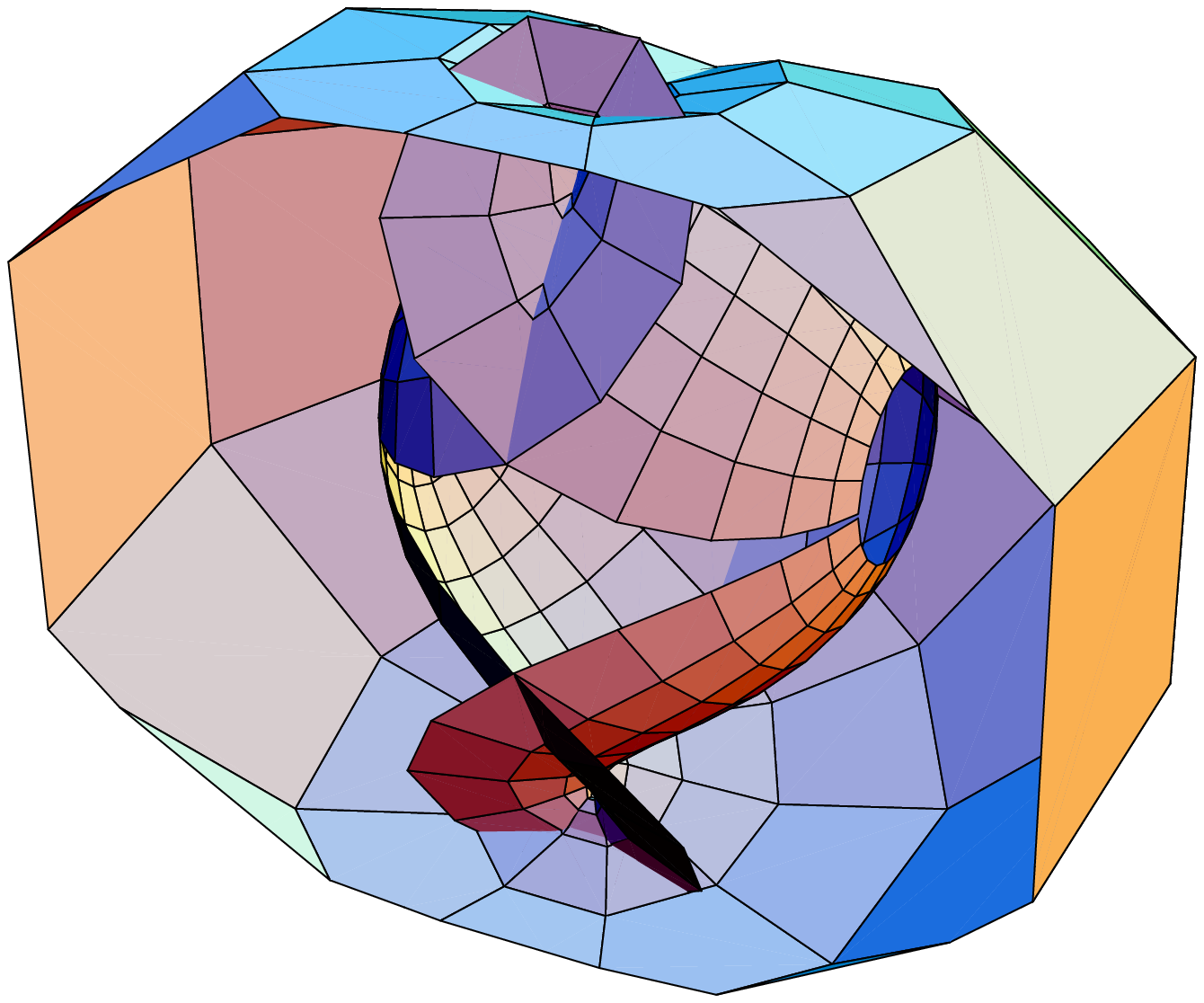}}
	{A D-transform of a spherical net: cmc-1 in $H^3$}\noindent}
\def\TPic{\Figure{ \hbox{
	\epsfxsize.30\hsize\epsfbox{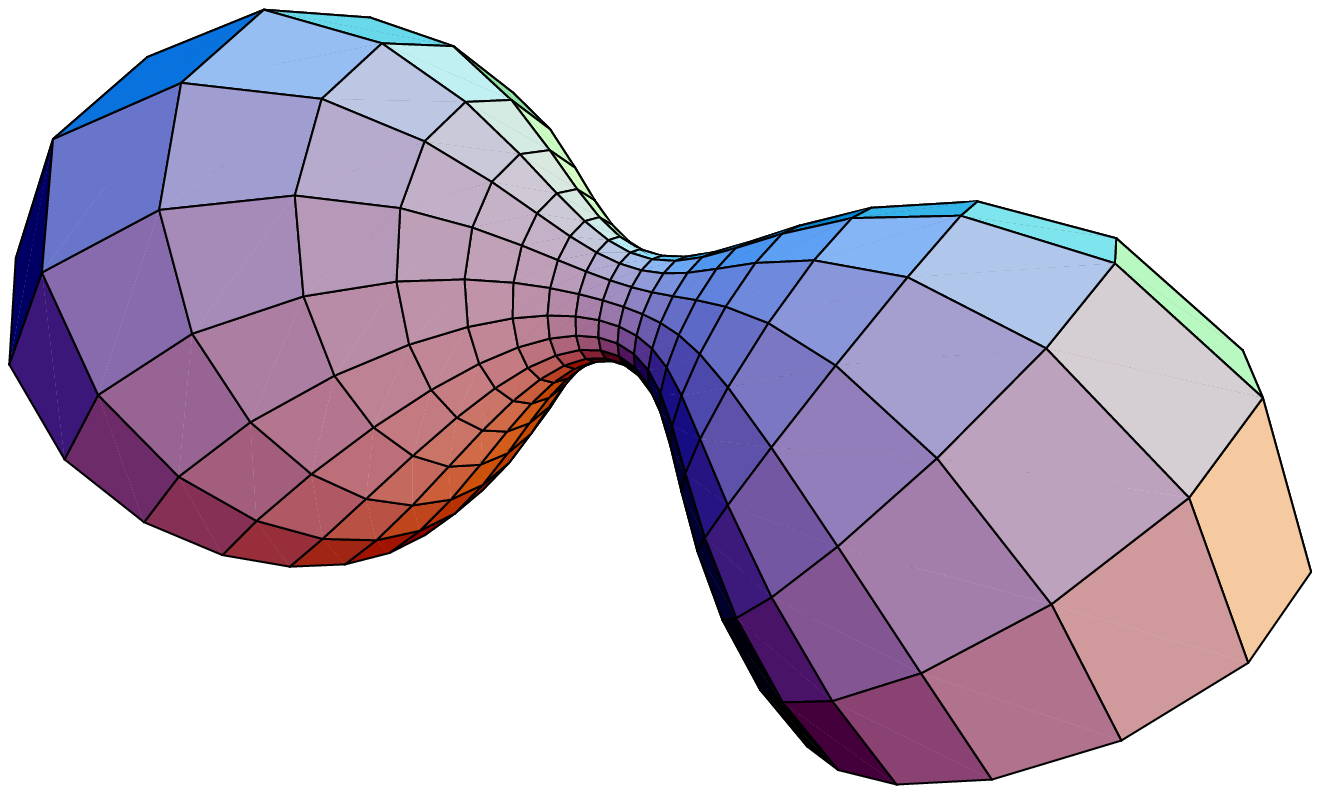} \vrule
	\epsfxsize.30\hsize\epsfbox{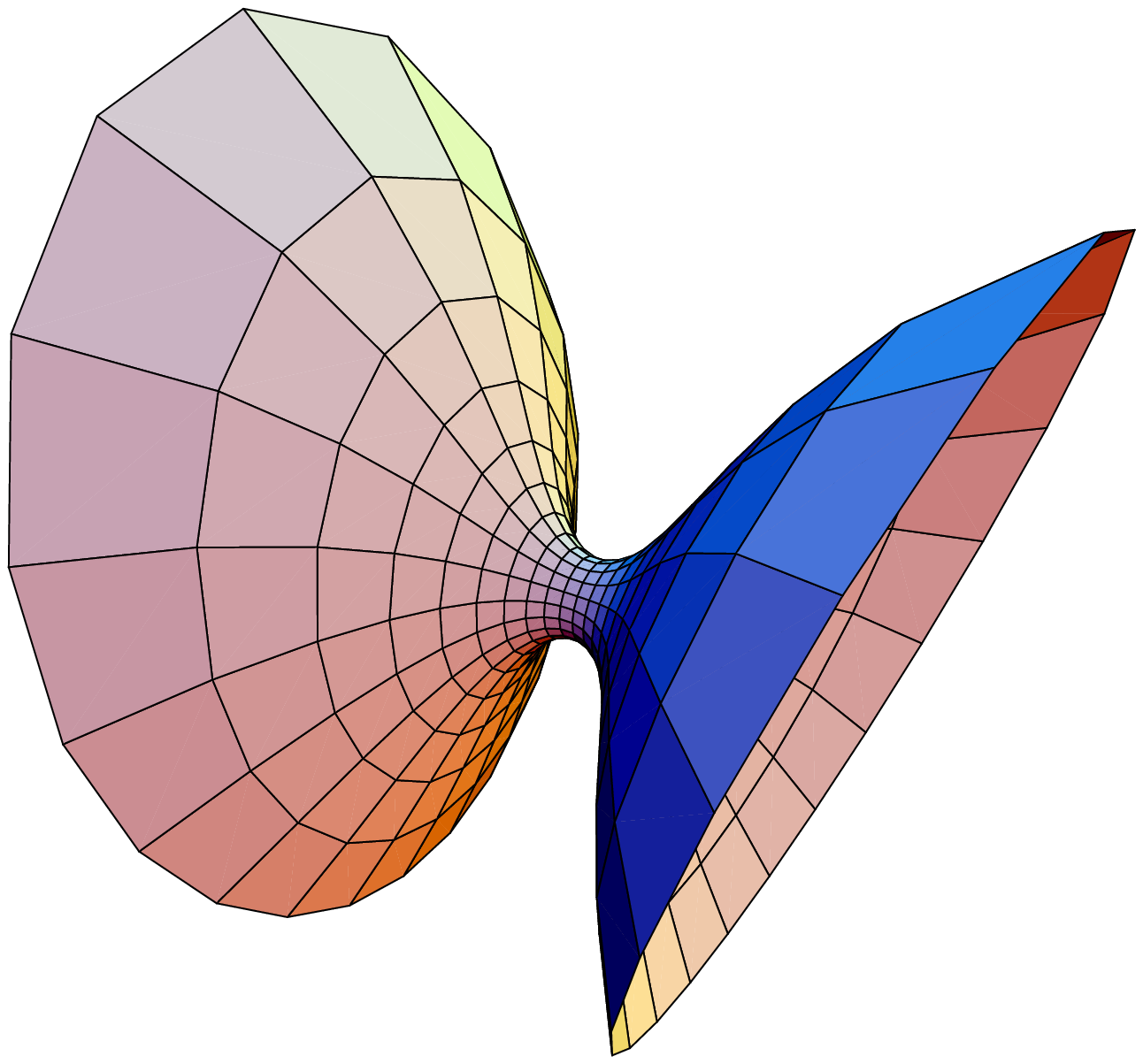} \vrule
	\epsfxsize.30\hsize\epsfbox{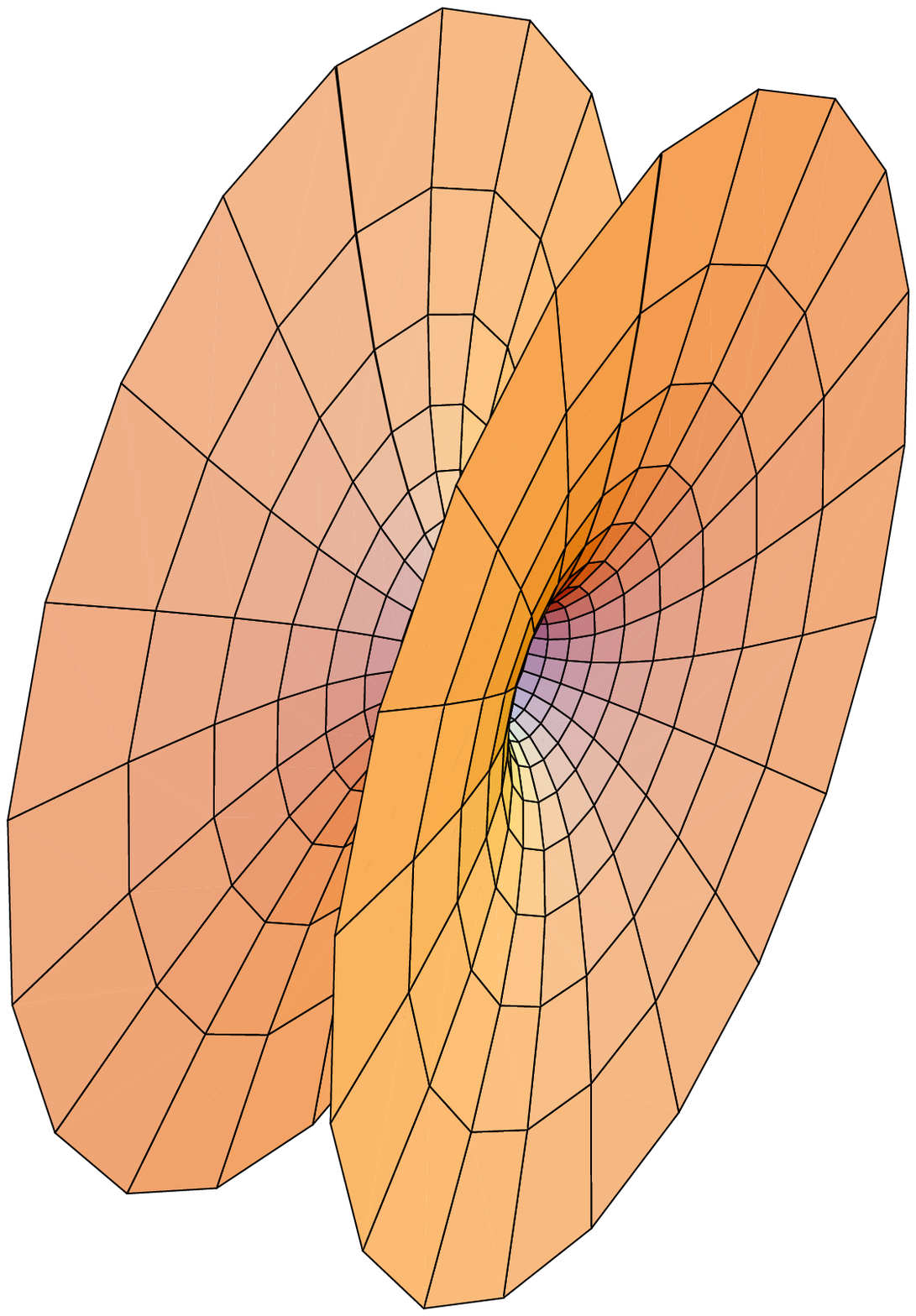}
	} \hrule \hbox{
	\epsfxsize.30\hsize\epsfbox{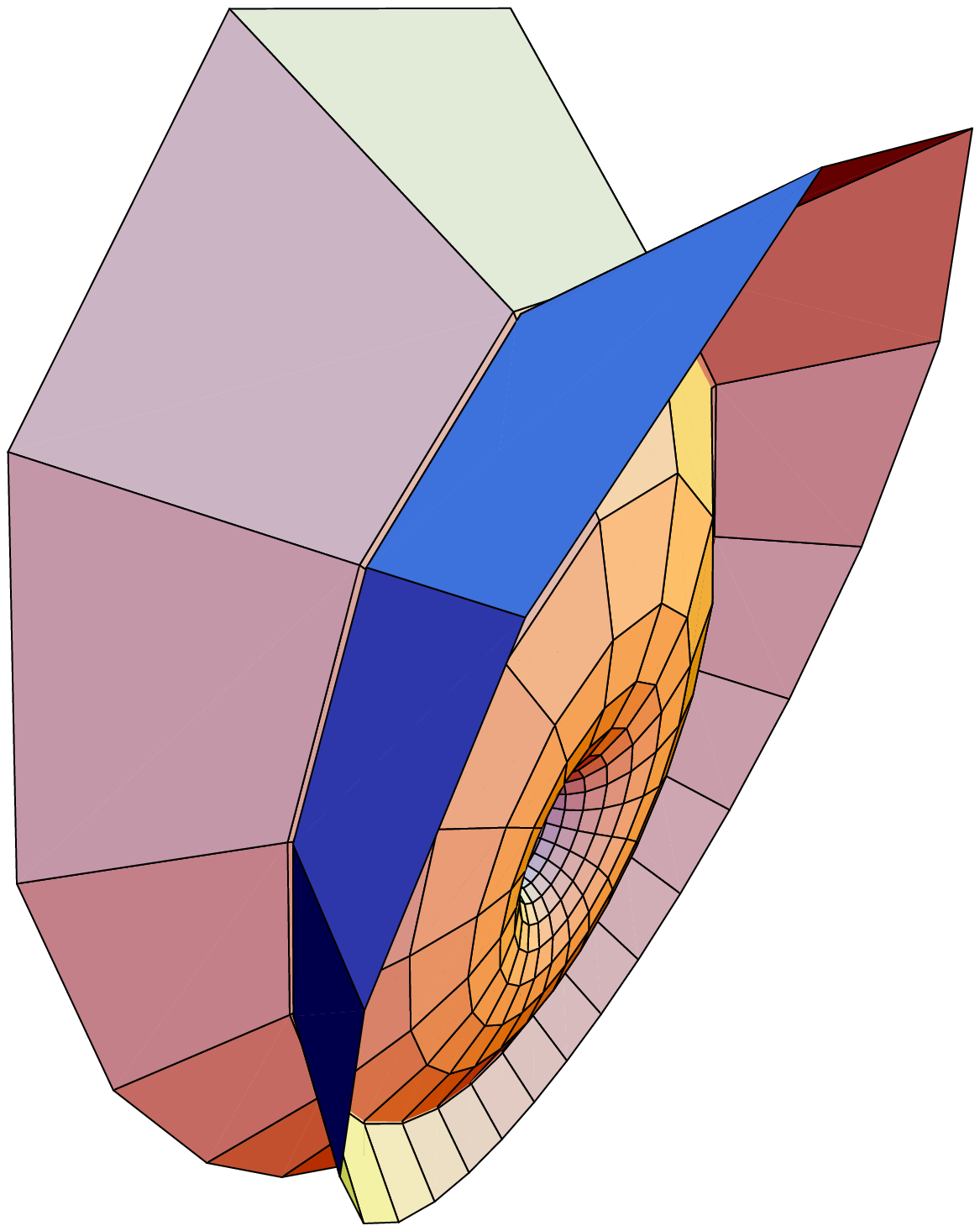} \vrule
	\epsfxsize.30\hsize\epsfbox{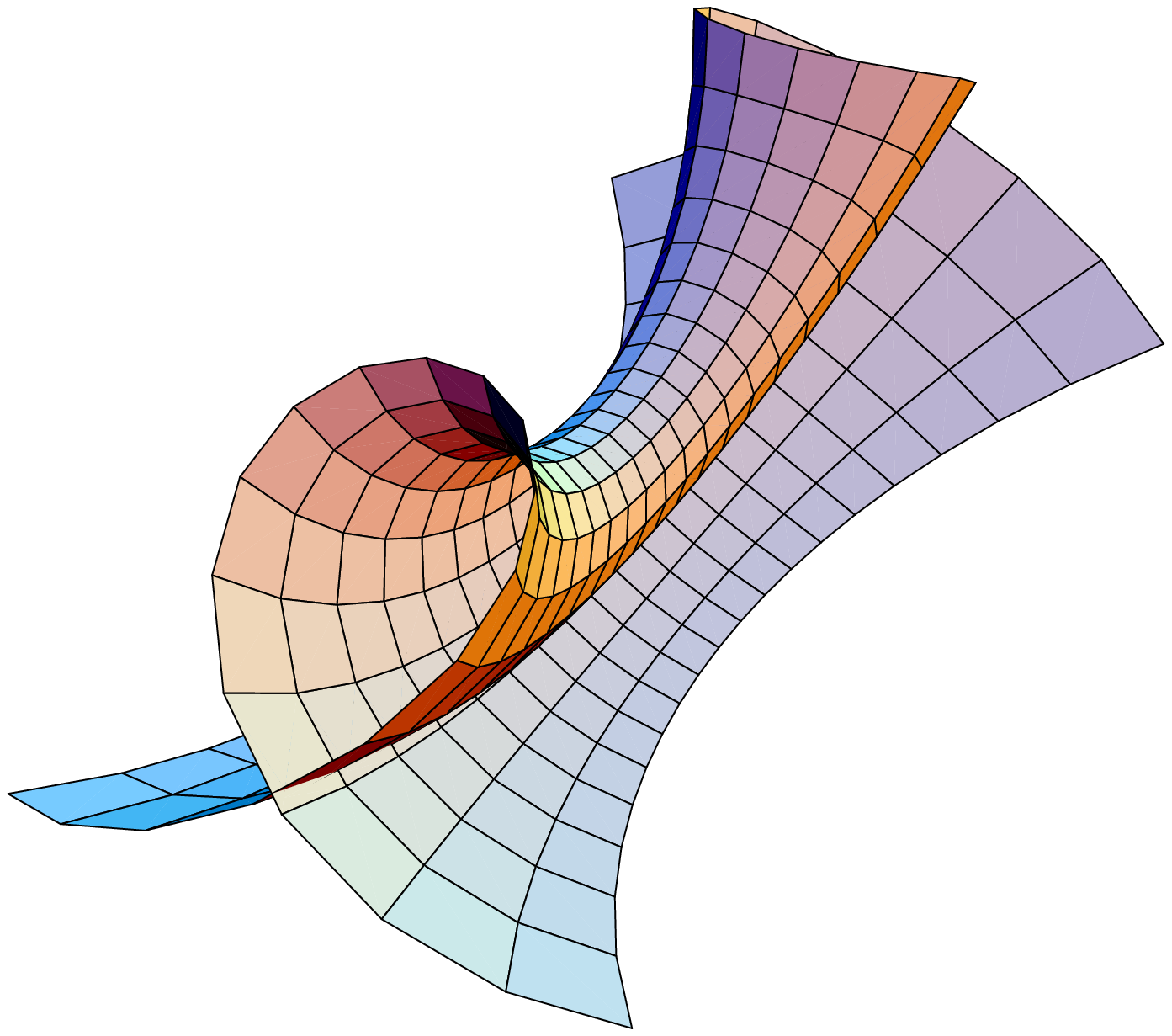} \vrule
	\epsfxsize.30\hsize\epsfbox{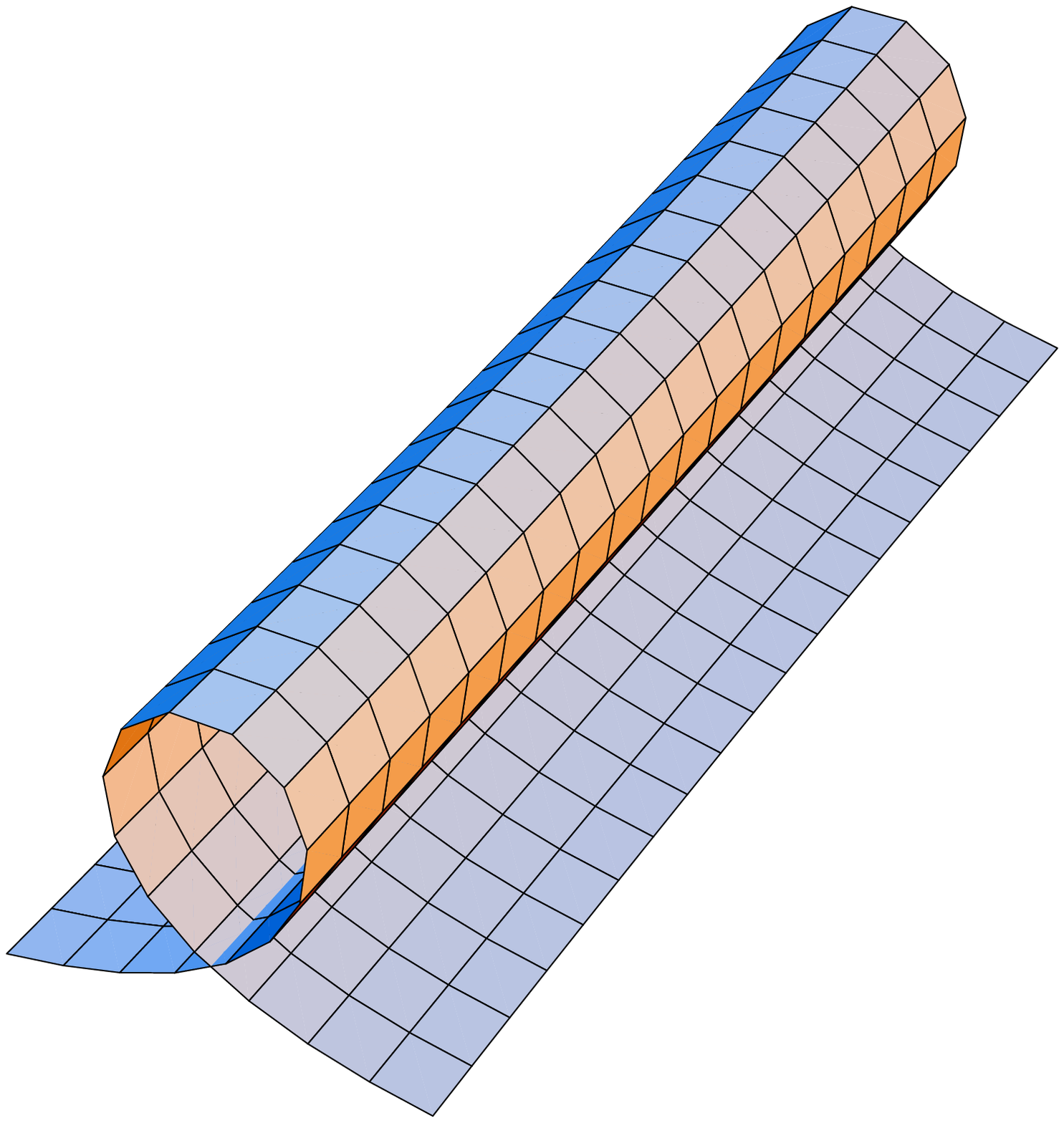}}
	}{T-transforms of the Catenoid: the Catenoid cousins}\noindent}


\BegRefs
\Ref\Asla H.~Aslaksen: {\it Quaternionic Determinants\/};
        Math.\ Intell.\ {\bf 18.3} (1996) 57-65
\Ref\Bianchi L.~Bianchi: {\it Complementi alle Ricerche sulle
       Superficie Isoterme\/}; Ann.\ Mat.\ pura appl.\ {\bf 12}
       (1905) 19--54
\Ref\Blaschke W.~Blaschke: {\it Vorlesungen \"uber Differentialgeometrie
        III\/}; Springer, Berlin 1929
\Ref\BoPi A.~Bobenko, U.~Pinkall: {\it Discrete Isothermic Surfaces\/};
        J.~reine angew.\ Math. {\bf 475} (1996) 187-208
\Ref\Bryant R.~Bryant: {\it Surfaces of mean curvature one in hyperbolic
        space\/}; Ast\'{e}risque {\bf 154-155} (1987) 321-347
\Ref\BJPP F.~Burstall, U.~Hertrich-Jeromin, F.~Pedit, U.~Pinkall:
        {\it Isothermic surfaces and Curved flats\/}; Math.~Z. {\bf 225}
        (1997) 199-209
\Ref\Calapso P.~Calapso: {\it Sulla superficie a linee di curvatura
	isoterme\/}; Rend. Circ. Mat. Palermo~{\bf 17} (1903) 275-286
\Ref\CalapsoT P.~Calapso: {\it Sulle trasformazioni delle superficie
	isoterme\/}; Annali Mat.~{\bf 24} (1915) 11-48
\Ref\Cartan \'{E}.~Cartan: {\it Les espaces \`{a} connexion conforme\/};
        Ann.\ Soc.\ Pol.\ Math.\ {\bf 2} (1923) 171-221
\Ref\Christoffel E.~Christoffel: {\it Ueber einige allgemeine
        Eigenschaften der Minimumsfl\"achen\/}; Crelle's J.\ {\bf 67}
        (1867) 218--228
\Ref\CGS J.~Cie\'{s}li\'{n}ski, P.~Goldstein, A.~Sym: {\it Isothermic
        Surfaces in $E^3$ as Soliton Surfaces\/}; Phys. Lett. A~{\bf 205}
        (1995) 37 - 43
\Ref\Darboux G.~Darboux: {\it Sur les surfaces isothermiques\/};
        Comptes Rendus {\bf 122} (1899) 1299--1305, 1483--1487, 1538
\Ref\FePe D.~Ferus, F.~Pedit: {\it Curved flats in Symmetric spaces\/};
        Manuscripta Math.\ {\bf 91} (1996) 445-454
\Ref\HHP U.~Hertrich-Jeromin, T.~Hoffmann, U.~Pinkall: {\it A discrete
        version of the Darboux transform for isothermic surfaces\/};
        to appear in A.~Bobenko, R.~Seiler, {\it Discrete integrable
        Geometry and Physics\/}, Oxford Univ.\ Press, Oxford 1997
\Ref\Suppl U.~Hertrich-Jeromin: {\it Supplement on Curved flats in the space
        of Point pairs and Isothermic surfaces: A Quaternionic calculus\/};
        Doc.\ Math.\ J.~DMV {\bf 2} (1997) 335-350
\Ref\JePe U.~Hertrich-Jeromin, F.~Pedit: {\it Remarks on the Darboux
        transform of isothermic surfaces\/};
        Doc.\ Math.\ J.~DMV {\bf 2} (1997) 313-333
\Ref\DiscreteT U.~Hertrich-Jeromin: {\it Transformations of discrete
	Isothermic surfaces and Discrete cmc-1 surfaces in hyperbolic
	space\/}; in preparation
\Ref\George G.~Kamberov: {\it Holomorphic quadratic differentials and
        Conformal immersions\/}; Preprint 1997
\Ref\MussoDef E.~Musso: {\it Deformazione di superficie nello spazio
	di M\"obius\/}; Rend.\ Istit.\ Mat.\ Univ.\ Trieste {\bf 27}
        (1995) 25-45
\Ref\Musso E.~Musso: {\it Isothermic surfaces in Euclidean space\/};
        L.~Cordero, E.~Garc\'{i}a-R\'{i}o (eds.), Proc.\ Workshop
        Recent Topics in Differential Geometry, July 1997, Public.\
        Depto.\ Geometr\'{i}a y Topolog\'{i}a, Univ.\ Santiago de
        Compostela {\bf 89} (1998) 219-235
\Ref\Joerg J.~Richter: {\it Conformal maps of a Riemann surface into the
	Space of Quaternions\/}; PhD thesis, TU Berlin, 1997
\Ref\Study E.~Study: {\it Ein Seitenst\"uck zur Theorie der linearen
        Transformationen einer komplexen Ver\-\"anderlichen, Teile I-IV\/};
        Math.~Z. {\bf 18} (1923) 55-86, 201-229 and {\bf 21} (1924)
        45-71, 174-194
\Ref\UmYa M.~Umehara, K.~Yamada: {\it A parametrization of the Weierstrass
        formulae and perturbation of complete minimal surfaces in $\R^3$
        into the hyperbolic 3-space\/}; J.~reine angew.\ Math.\ {\bf 432}
        (1992) 93-116
\Ref\UmYaDual M.~Umehara, K.~Yamada: {\it A Duality on cmc-1 surfaces in
	hyperbolic space and a hyperbolic analogue of the Ossermann inequality\/};
	Tsukuba J.\ Math.\ {\bf 21} (1997) 229-237
\Ref\Wilker J.~Wilker: {\it The Quaternion formalism for M\"obius groups
        in four or fewer dimensions\/}; Lin.\ Alg.\ Appl.\ {\bf 190} (1993)
        99-136
\EndRefs


\Section{Introduction}
This paper aims to establish a relation between two fields of current interest.

On one side, isothermic surfaces and their various transformations gained new
interest when an integrable system approach was developed \CGS: ``Darboux
pairs'' of isothermic surfaces correspond to ``curved flats'' in a suitable
Grassmannian \BJPP\ (cf.\Suppl), a particularly simple type of integrable
system. The ``spectral parameter'' occuring in the loop group description
for curved flats \FePe\ gives rise to a ``spectral transformation'' for
isothermic surfaces. Classically, this transformation was independently
discovered by P.~Calapso \Calapso, \CalapsoT\ who considered it the most
fundamental transformation, and by L.~Bianchi \Bianchi\ who introduced it
as the ``$T$-transformation''. Later, isothermic surfaces were identified as
the only surfaces in conformal geometry allowing a second order deformation
\Cartan\ --- that turns out to be exactly the $T$-transformation \MussoDef.
Another central role is played by the ``Darboux transformation'' for isothermic
surfaces, first introduced by G.~Darboux \Darboux: a corresponding
permutability theorem and a quaternionic calculus \Suppl\ are the
keys to the theory of ``discrete isothermic nets'' \BoPi, and their
transformations \HHP, \DiscreteT.

On the other side, surfaces of constant mean curvature 1 in hyperbolic space 
became a field of interest, complementing research on minimal surfaces in
Euclidean space, after R.~Bryant \Bryant\ developed a Weierstrass type
representation for them. Subsequently, M.~Umehara and K.~Yamada described
these surfaces (and the corresponding Weierstrass representation) as
a deformation of minimal surfaces in Euclidean space (and their Weierstrass
representation) \UmYa, and introduced a notion of duality \UmYaDual\ for
cmc-1 surfaces in $H^3$.

In section 3 of the present paper, we describe the various transformations of
isothermic surfaces in terms of a quaternionic calculus that we shortly present
in section 2 (for more information see \Suppl). Here, we emphasize the
geometric context of the transformation, e.g.\ the Christoffel transformation
depends on the Euclidean geometry of the ambient space while the Goursat
transformation uses the interplay of the ambient space's Euclidean geometry
and the conformal invariance of the notion of isothermic surfaces; the
Darboux transformation is a transformation for surfaces in M\"obius space
while the $T$-transformation is only well defined for M\"obius equivalence
classes of isothermic surfaces. Then, we discuss the interrelations of these
transformations in terms of several ``permutability theorems'' --- some of
the proofs become very compact in the quaternionic calculus.

In section 4, we interpret certain results on cmc-1 surfaces in
hyperbolic space in terms of transformations of isothermic surfaces:
the Umehara-Yamada deformation of minimal into cmc-1 surfaces is identified
as the $T$-transformation and a unified version of the Weierstrass type
representations for both surface classes is given in terms of quaternions.
In this context, we give a M\"obius geometric characterization for cmc-1
and minimal surfaces --- in fact, the Umehara-Yamada perturbation families
of cmc-1 surfaces appear as spectral families of isothermic surfaces. In terms
of M\"obius geometry, this means that the induced metric of the central sphere
congruence has constant Gauss curvature 1.
Finally, we introduce a second Weierstrass type representation for cmc-1
surfaces in $H^3$ in terms of the Darboux transformation, and discuss its
relation with the notion of duality for cmc-1 surfaces.

As all the transformations under discussion require an integration they generally
do not preserve periods, and may develop singularities, too. In this respect,
our results have to be understood as local results. However, we formulate
most of them in an invariant way in order to make the theory we present
more accessible for further research on the global geometry of isothermic
surfaces and their transformations. For example, it should be an interesting
task to examine to what extent the Darboux and $T$-transformations respect
extrinsic symmetries of surfaces, or how they behave in the neighbourhood
of an isolated umbilic.

As mentioned above, many of the results of this paper can be formulated
in a very elegant way using

\vfill\eject 

\Section{The Quaternionic formalism}
for M\"obius differential geometry, as introduced in \Suppl\
(cf.\Study,\Wilker) --- some results even {\it depend} on that approach.
On the other hand, the classical approach of considering M\"obius
geometry as a subgeometry of {\it real} projective geometry is the
better known formalism and allows a more direct treatment of certain
aspects in local differential geometry.
In this section, we intend to discuss the relation between the two
approaches briefly, focussing on 3-dimensional geometry\Footnote{For
a more comprehensive treatment --- particularly for the 4-dimensional
setting that also appears natural when using quaternions --- the reader
is referred to \Suppl.}.

Classically, the conformal $n$-sphere $S^n$ is embedded into the real
projective $(n+1)$-space $\RP^{n+1}$. Then, the M\"obius transformations
of $S^n$ become the projective transformations of $\RP^{n+1}$ that fix
$S^n\subset\RP^{n+1}$ as an absolute quadric.
Up to scaling, there is a unique Lorentz scalar product on the space
$\R^{n+2}_1$ of homogeneous coordinates such that the points of $S^n$
are exactly the isotropic lines of that Lorentz product.
Via orthogonality (polarity), hyperspheres can be identified with
spacelike lines in $\R^{n+2}_1$, and the M\"obius transformations
appear as pseudo orthogonal transformations (cf.\Blaschke).

On the other hand, the orientation preserving M\"obius transformations
of the conformal 4-sphere $\H\cup\{\infty\}$ can be identified with
the projective transformations of the quaternionic projective line
$\HP^1$ --- just as the (orientation preserving) M\"obius transformations
of the conformal 2-sphere $\C\cup\{\infty\}$ can be identified with the
projective transformations of $\CP^1$.
Thus, in homogeneous coordinates\Footnote{We consider $\H^2$
as a {\it right} vector space over the quaternions.} $v\in\H^2$ for
points $v\H\in\HP^1\cong S^4$, the orientation preserving M\"obius
transformations appear as linear transformations $v\mapsto Av$ with
$A\in Gl(2,\H)$. Or, in affine coordinates $x\simeq(x,1)^t$, the
M\"obius transformations become fractional linear transformations
$x\mapsto(ax+b)(cx+d)^{-1}$ on $\H\cong\R^4$.

The relation between the two models can be established using the space
of {\it quaternionic hermitian forms} $s:\H^2\times\H^2\to\H$ on $\H^2$,
$$s(u,v\lambda+w\mu)=s(u,v)\lambda+s(u,w)\mu\hskip1em{\rm and}
	\hskip1em s(w,v)=\overline{s(v,w)}.$$
Fixing a basis $(e_1,e_2)$ of the quaternionic plane $\H^2$, any
quaternionic hermitian form is determined by its values $s(e_1,e_1),
s(e_2,e_2)\in\R$ and $s(e_1,e_2)\in\H$. Thus, the space of quaternionic
hermitian forms is a 6-dimensional real vector space. Endowed with
the determinant as a quadratic form\Footnote{Indeed, the Lorentz
product introduced this way is only well defined up to a scaling:
a change of basis $(e_1,e_2)\to(\tilde{e}_1,\tilde{e}_2)$ results
in a rescaling of $\langle s,s\rangle\to\lambda^2\langle s,s\rangle$
(cf.\Suppl).  This is the same scaling ambiguity that arises in the
classical construction, above.},
$$\langle s,s\rangle := -\det\left(\matrix{s(e_1,e_1)&s(e_1,e_2)\cr
	s(e_2,e_1)&s(e_2,e_2)\cr}\right) 
	= |s(e_1,e_2)|^2-s(e_1,e_1)\cdot s(e_2,e_2),$$
this space becomes the 6-dimensional Minkowski space $\R^6_1$ of the
classical model for 4-dimensional M\"obius geometry: space- and
lightlike lines of quaternionic hermitian forms can be identified
with their null cones in $\H^2$ to encode hyperspheres $s\subset\HP^1$
and points $p\in\HP^1$ (cf.\Suppl,\Study).
Moreover, using the Study determinant ${\cal D}$ (cf.\Asla) for
quaternionic $2\times2$-matrices, the special linear group
$$Sl(2,\H):=\{M\in M(2\times2,\H)\,|\,{\cal D}M=1\}$$
acts by isometries on $\R^6_1$ via $(M,s)\to M\cdot s:=s(M^{-1}.,M^{-1}.)$.
In fact, $Sl(2,\H)$ is the universal cover of the identity component of
the M\"obius group.

As hyperspheres $s\subset\HP^1$ are encoded as spacelike lines of
quaternionic hermitian forms, the M\"obius group of the conformal
3-sphere can be obtained as the subgroup of $Sl(2,\H)$ fixing a given
quaternionic hermitian form, say\Footnote{Fixing one of the two length
1 representatives of the spacelike line amounts to fixing an orientation
on the conformal 3-sphere $S^3$.}
$$S^3\cong\left(\matrix{0&1\cr1&0\cr}\right).$$
In affine coordinates, this choice corresponds to the imaginary
quaternions as Euclidean 3-space, ${\rm Im}\H\cong\R^3$: the null
cone of the above form is
$$S^3 = \{\left(\matrix{h\cr1\cr}\right)\cdot\H\,|\,h\in{\rm Im}\H\}
      \cup\{\left(\matrix{1\cr0\cr}\right)\cdot\H\}\subset\HP^1.$$
In this setting, the Lie algebra of the M\"obius group of $S^3$ becomes
$${\fk m} = \{\left(\matrix{r+a&b\cr c&-r+a\cr}\right)\,|\,
        r\in\R; a,b,c\in{\rm Im}\H\} \subset {\fk sl}(2,\H).$$
Space- and lightlike quaternionic hermitian forms that are perpendicular
to the above fixed form representing $S^3$ can now be identified with
2-spheres and points, respectively, in the conformal 3-sphere $S^3$.

Given a surface $f:M^2\to S^3$, a congruence of 2-spheres
$s:M^2\to\R^6_1$ in $S^3$ is said to be {\it enveloped} by $f$ if
each point $f(p)$ lies on the corresponding sphere $s(p)$, and if
the tangent planes $d_pf(T_pM)=T_{f(p)}s(p)$ coincide.
As discussed, incidence is encoded by $s(p)(f(p),f(p))=0$ in terms
of homogeneous coordinates\Footnote{We will use the same letter $f$
for the immersion into $S^3\subset\HP^1$ and for its homogeneous
coordinates into $\H^2$.} $f:M^2\to\H^2$. A simple calculation (cf.\Suppl)
shows that

\Lemma{}{An immersion $f:M^2\to S^3$ envelopes a congruence $s:M^2\to\R^6_1$
of 2-spheres in $S^3$ if and only if
$$s(f,f)\equiv0\hskip1em{\rm and}\hskip1em
        s(df,f)+s(f,df)\equiv0.\Eqno\Label\Envelope$$}

\noindent
Fixing a basis $(e_1,e_2)$ of $\H^2$, a map $F:M^2\to{\rm Sl}(2,\H)$ will
be called an {\it adapted frame} of an immersion $f:M^2\to S^3$ if $f=Fe_1$.
And, $F$ is an adapted frame for a {\it surface pair} $f,\hat{f}:M^2\to S^3$
if $f=Fe_1$ and $\hat{f}=Fe_2$.
Note, that for a surface pair $(f,\hat{f})\simeq F$ the fact that $F$ is
invertible is equivalent to $f(p)\neq\hat{f}(p)$ at all points $p\in M^2$.
Denoting $s_f,s_{\hat{f}}\in\R^6_1$ two generators of the lightlike lines
of quaternionic hermitian forms corresponding to $f$ and $\hat{f}$,
respectively, $s_f(f,f)=0$ together with $|s_f|^2=0$ implies that also
$s_f(f,\hat{f})=0$. Thus, choosing a suitable scaling for $s_f$ and
$s_{\hat{f}}$, we have
$$F^{-1}\cdot s_f=\left(\matrix{0&0\cr0&1\cr}\right)\hskip1em{\rm and}
        \hskip1em F^{-1}\cdot s_{\hat{f}}=\left(\matrix{1&0\cr0&0\cr}\right).
        \Eqno\Label\ClassicalFrameF$$
If $s:M^2\to\R^6_1$ is a sphere congruence in $S^3$ containing the points
of both surfaces of a surface pair $(f,\hat{f})$, then the first condition
in \Envelope\ yields $$F^{-1}\cdot s=\left(\matrix{0&c\cr-c&0\cr}\right)$$
with a function $c:M^2\to{\rm Im}\H$ into the imaginary quaternions.
Given three such sphere congruences, say $s_i,s_j,s_k:M^2\to\R^6_1$
defined by
$$F^{-1}\cdot s_i=\left(\matrix{0&i\cr-i&0\cr}\right),\hskip1em
        F^{-1}\cdot s_j=\left(\matrix{0&j\cr-j&0\cr}\right),\hskip1em
        F^{-1}\cdot s_k=\left(\matrix{0&k\cr-k&0\cr}\right),
        \Eqno\Label\ClassicalFrameS$$
a classical frame for M\"obius geometry in $S^3$ (cf.\Blaschke) is obtained:
it consists of three congruences of orthogonally intersecting spheres $s_i$, $s_j$,
and $s_k$, and the maps $s_f\cong f$ and $s_{\hat{f}}\cong\hat{f}$ describing
their intersection points.

Since $d(F^{-1}\cdot s_c)=0$ for all the above defined frame vectors $s_c$,
one calculates $ds_c=-s_c(.,\Phi.)-s_c(\Phi.,.)$ with the connection form
$\Phi=F^{-1}dF$.
Using this equation allows to write the structure equations of the classical
frame in terms of the components of $\Phi$ --- thus obtaining a ``translation
table'' between the classical and the quaternionic setting (in the codimension
1 case):
$$\Phi={\petit\matrix{\left(\matrix{
	\langle f,d\hat{f}\rangle
	+{\langle ds_j,s_k\rangle\over2} i
	+{\langle ds_k,s_i\rangle\over2} j
	+{\langle ds_i,s_j\rangle\over2} k
	&
	 \langle ds_i,\hat{f}\rangle i
	+\langle ds_j,\hat{f}\rangle j
	+\langle ds_k,\hat{f}\rangle k
	\cr
	 \langle df,s_i\rangle i
	+\langle df,s_j\rangle j
	+\langle df,s_k\rangle k
	&
	\langle df,\hat{f}\rangle
	+{\langle ds_j,s_k\rangle\over2} i
	+{\langle ds_k,s_i\rangle\over2} j
	+{\langle ds_i,s_j\rangle\over2} k
	\cr}\right)}}.$$
For example, the surface $f=Fe_1$, or $\hat{f}=Fe_2$, will envelope the
sphere congruence $s_i$ if and only if the $i$-component $\langle df,s_i
\rangle=0$, or $\langle d\hat{f},s_i\rangle=0$ respectively, in the
corresponding off-diagonal form of $\Phi$ vanishes (cf.\Envelope).

At this point, we are prepared to collect (and reformulate) some results on

\Section{Transformations of Isothermic surfaces}
Classically, an immersed surface $f:M^2\to\R^3\cong{\rm Im}\H$ is called
``isothermic'' if it carries conformal curvature line coordinates
$(x,y):M^2\to\R^2$, i.e.\ $$\matrix{
	\langle df,df\rangle &=& e^{2u}(dx^2+dy^2),\hfill\cr
	\langle df,dn\rangle &=& -e^{2u}(k_1dx^2+k_2dy^2),\cr}$$
where $n:M^2\to S^2$ is the normal field of $f$ and $u:M^2\to\R$ is a
suitable real valued function.
This definition obviously runs into problems when umbilics are present
--- and, since we are planning to consider sphere pieces as isothermic surfaces,
we look for a more appropriate definition.
On the tangent bundle of an immersed surface $f:M^2\to{\rm Im}\H$,
(quaternionic) left multiplication with the unit normal field $n:M^2\to S^2$
yields a $90^0$-rotation; the pull back of this structure provides a complex
structure $J$ on $M^2$: $df(Jx)=n\cdot df(x)$. Clearly, the immersion $f$
is conformal with respect to this complex structure.
It is well known that, away from umbilics, the existence of conformal
curvature line coordinates is equivalent to the fact that the Hopf
differential $\eta$ of a conformal immersion $f$ is a real multiple
of a holomorphic quadratic differential $q$ on $M^2$: $\eta=\varrho q$
with $\varrho:M^2\to\R$.
Thus, we give the definition of an isothermic immersion of
a Riemann surface with ``prescribed'' principal curvature lines:

\Definition{~(isothermic immersion)}{A Riemann surface $M^2$ equipped
with a holomorphic quadratic differential\Footnote{In order to establish
a theory of ``globally isothermic surfaces'' it seems too restrictive to
ask $q$ to be holomorphic --- for example, with this assumption the
ellipsoid is no longer isothermic \George. In order to give a
satisfactory global definition, the consequences of various possible
assumptions on $q$ have still to be worked out.} $q$ is called a
polarized surface; \par\noindent
a conformal immersion $f:(M^2,q)\to{\rm Im}\H$ of a polarized surface
is called isothermic if its Hopf differential $\eta=\varrho q$ is
a real multiple, $\varrho:M^2\to\R$, of $q$.}

\SubSection{The Christoffel transformation}
If two non homothetic immersions $f,f^{\ast}:M^2\to{\rm Im}\H$
induce the same conformal structure and have parallel normal fields
then \Christoffel\ they are either associated minimal surfaces
($n^{\ast}=n$), or $n^{\ast}=-n$ and both surfaces are isothermic.
The latter case can be characterized (cf.\Suppl) by the equation
$$df\wedge df^{\ast}=0.\Eqno\Label\CEqn$$
On the other hand, given an isothermic surface in terms of isothermic
curvature line coordinates, $f:(\C,dz^2)\to{\rm Im}\H$,
$df^{\ast}:=f_x^{-1}dx-f_y^{-1}dy$ defines, at least locally, a second
surface satisfying $df\wedge df^{\ast}=0$.
Obviously, the surface $f^{\ast}$ is only determined up to homothety
(and translation) by \CEqn. Canonically trivializing the complex line
bundle ${\rm span}\{1,n\}\cong M^2\times\C$ and using the above scaling
of $f^{\ast}$, one obtains the polarization $q$ back:
$$df\cdot df^{\ast}=(dx^2-dy^2)+2n\,dxdy\simeq dz^2.$$
In fact, if $\tau:TM^2\to{\rm Im}\H$ is a 1-form with values
in the imaginary quaternions then, $df\wedge\tau=0$ if and only if
$\tau(Jx)+n\tau(x)=\tau(Jx)-\tau(x)n=0$ if and only if
$(df\cdot\tau)((a+bJ)x)=(a+bn)^2(df\cdot\tau)(x)$ and $n\tau(x)+\tau(x)n=0$
for any $x\in TM$:

\Lemma{}{Let $f:M^2\to{\rm Im}\H$ be a conformal immersion of a Riemann
surface and $\tau:TM\to{\rm Im}\H$ a 1-form with values in the imaginary
quaternions. Then $df\cdot\tau:TM\to{\rm span}\{1,n\}$ defines a
quadratic differential if and only if $df\wedge\tau=0$.}

\noindent
This lemma allows us (cf.\George) to define {\it the} Christoffel
transform of a polarized isothermic surface (up to translation),
i.e.\ to fix the scaling of $f^{\ast}$:

\Definition{~(Christoffel transformation)}{Let $f:(M^2,q)\to{\rm Im}\H$
be an isothermic immersion of a polarized Riemann surface.
Then $\Ct{f}:M^2\to{\rm Im}\H$ is called Christoffel transform of $f$ if
$df\cdot d\Ct{f}\simeq q$.}

\noindent
The Christoffel transform $\Ct{f}$ of an isothermic immersion $f$ is
unique up to translations and it is isothermic with respect to the
polarization $\bar{q}=d\Ct{f}\cdot df$.
The Christoffel transform of $\Ct{f}$ is the original surface, $f=C^2f$.

As an example, we consider a minimal surface $f:M^2\to{\rm Im}\H$ in
Euclidean 3-space:
since $df\wedge dn=-Hdf\wedge df$, minimal surfaces are exactly those
isothermic surfaces whose Christoffel transforms are totally umbilic
(cf.\JePe). As all the surfaces $f_t$ in the associated family,
$df_t=(\cos t+\sin t\cdot n)df$, have the same Gauss map the additional
information provided by the (holomorphic) Hopf differential $df\cdot dn$
as polarization is required to reconstruct the original surface from its
Christoffel transform (Gauss map).

\SubSection{The Goursat transformation}
is classically defined for minimal surfaces, only: the action of the
complex orthogonal group $O(3,\C)$ on the holomorphic null curve
describing an (associated family of) minimal surfaces provides a
3-parameter family of nontrivial transformations.
For the Gauss map (Christoffel transform) of the minimal surface,
these correspond to nontrivial M\"obius transformations.
Obviously, this type of transformation can be generalized to isothermic
surfaces ---
making use of the fact that the notion of isothermic surface is conformally
invariant while the Christoffel transformation depends on Euclidean geometry
of the ambient space:

\Definition{~(Goursat transformation)}{Let $f:(M^2,q)\to{\rm Im}\H$ be
isothermic, and let $M:{\rm Im}\H\cup\{\infty\}\to{\rm Im}\H\cup\{\infty\}$
be a M\"obius transformation.
Then an immersion $\Gt{f}:=CMCf$, i.e.\ the Christoffel transform of the
M\"obius transformed Christoffel transform of $f$, is called a Goursat
transform of $f$.}

\noindent
The Goursat transformations form a group acting on isothermic immersions
of a polarized Riemann surface since $C^2=id$.
If the M\"obius transformation $M$ is a similarity, then the corresponding 
Goursat transform $\Gt{f}$ of an isothermic immersion is clearly similar to
the original immersion $f$.
If the M\"obius transformation $M$ is ``essential'', however, $\Gt{f}$ will
generally be not M\"obius equivalent to the original surface. In this way,
the Goursat transforms of an isothermic immersion $f$ provide a 3-parameter
family of non M\"obius equivalent isothermic immersions.
Thus, we may restrict the attention to essential M\"obius transformations
--- that are (up to similarity) of the form $x\to(x-m)^{-1}$ with
$m\in{\rm Im}\H$:

\Lemma{}{Let $x\to M(x)=(x-m)^{-1}$, $m\in{\rm Im}\H$, be an essential
M\"obius transformation, and $f:(M^2,q)\to{\rm Im}\H$ an isothermic
immersion.
Then \Suppl, $$d\Gt{f}=-(\Ct{f}-m)\cdot df\cdot(\Ct{f}-m).\Eqno\Label\GEqn$$}

\noindent
Consequently, the Goursat transformation appears as a special type of
spin transformation \Joerg.
As an example, we obtain the classical Weierstrass representation of
minimal surfaces: given a meromorphic function $g:M^2\to\C$ and a holomorphic
differential $\omega:TM^2\to\C$, the two maps\Footnote{Here, we do not worry
about periodicity of $\int\!\omega$, we {\it assume} the integral to exist.}
$\int\!\omega j,-jg:M^2\to\C j$ form a (trivial) Christoffel pair with respect
to the polarization $q=\omega\cdot dg$.
With the stereographic projection $\C j\ni x\to-i-2(i+x)^{-1}\in S^2$,
\GEqn\ yields the representation
$$df={1\over2}(i-jg)\omega j(i-jg)\Eqno\Label\WeierstrassRepresentation$$
of a minimal surface $f:M^2\to{\rm Im}\H$ with Weierstrass data $(g,\omega)$
--- cf.\Suppl.

\SubSection{The Darboux transformation}
A sphere congruence enveloped by two surfaces induces a point-to-point
correspondence between the two surfaces. If, under this correspondence,
the induced \DPic metrics on both surfaces are conformally equivalent and their
curvature lines do correspond, then \Darboux\ either the two surfaces are
M\"obius equivalent or, both surfaces are isothermic.
In terms of parametrizations $f,\hat{f}:M^2\to{\rm Im}\H$, the latter
situation can be characterized by the equation
$$df\wedge(f-\hat{f})^{-1}d\hat{f}=0.\Eqno\Label\DEqn$$
Comparison with \CEqn\ leads (cf.\JePe) to the following

\Lemma{}{Let $f:M^2\to{\rm Im}\H$ be an immersion. If there is a second
immersion $\hat{f}:M^2\to{\rm Im}\H$ satisfying \DEqn\ then $f$ is isothermic
with $q=df\cdot d\Ct{f}$ where
$$\lambda d\Ct{f}=(\hat{f}-f)^{-1}d\hat{f}(\hat{f}-f)^{-1}\Eqno\Label\CDEqn$$
with some $\lambda\in\R\setminus\{0\}$ defines the Christoffel transform of
$f$.}

\noindent
Rewriting \CDEqn\ yields the Riccati equation $d(\hat{f}-f)=(\hat{f}-f)\lambda
d\Ct{f}(\hat{f}-f)-df$ for the difference vector field of $f$ and $\hat{f}$
(cf.\JePe). Given an isothermic immersion $f:(M^2,q)\to{\rm Im}\H$, this
Riccati equation is completely integrable, i.e.\ fixing the parameter $\lambda$
and an initial condition there exists (at least locally) a solution $\hat{f}$.
This allows us to define the Darboux transforms $D_{\lambda}$ of an isothermic
immersion:

\Definition{~(Darboux transformation)}{Let $f:(M^2,q)\to{\rm Im}\H$ be an
isothermic immersion of a polarized Riemann surface. Then $\Dt{\lambda}{f}:M^2\to{\rm
Im}\H$ is called a $\lambda$-Darboux transform of $f$ if
$$d(\Dt{\lambda}{f}-f) =(\Dt{\lambda}{f}-f)\lambda d\Ct{f}(\Dt{\lambda}{f}-f)-df.\Eqno\Label\Riccati$$}

\noindent
The linear version $0=dv+\Phi_{\lambda}v$, $v:M^2\to\H^2$, with $\Phi_{\lambda}
=\left(\matrix{0&\lambda d\Ct{f}\cr df&0\cr}\right)$ of \Riccati\ is exactly
Darboux's linear system \Darboux: given a solution $v=(v_1,v_2)^t$, the
corresponding Darboux transform\Footnote{Choosing an initial condition
$v(p_0)\in S^3$ the solution will stay in $S^3$, thus $\Dt{\lambda}{f}:M^2\to{\rm Im}\H$.}
of $f$ is given by $\Dt{\lambda}{f}=f+v_2v_1^{-1}$.
In terms of homogeneous coordinates this reads $\Dt{\lambda}{f}=F_0v$ with the canonical
Euclidean frame $F_0=\left(\matrix{f&1\cr1&0\cr}\right)$. This description
reflects the invariance of the notion of a Darboux transform under the action
of the M\"obius group: if $M\in Sl(2,\H)$ is a M\"obius transformation then
$\Dt{\lambda}{Mf}=MF_0v=M\Dt{\lambda}{f}$.

\SubSection{The T-transformation}
Obviously, the above 1-form $\Phi_{\lambda}:TM^2\to{\fk sl}(2,\H)$ satisfies
the Maurer-Cartan equation $d\Phi_{\lambda}+\Phi_{\lambda}\wedge\Phi_{\lambda}
=0$. Consequently, it can (locally) be integrated to a frame $F_{\lambda}:
M^2\to Sl(2,\H)$, $dF_{\lambda}=F_{\lambda}\Phi_{\lambda}$.

\Definition{~(T-transformation)}{Let $f:(M^2,q)\to{\rm Im}\H$ be an isothermic
immersion of a polarized Riemann surface, and let $F_{\lambda}:M^2\to Sl(2,\H)$
be a frame with connection form $F_{\lambda}^{-1}dF_{\lambda}=\Phi_{\lambda}$,
$$\Phi_{\lambda}=\left(\matrix{0&\lambda d\Ct{f}\cr df&0\cr}\right).
	\Eqno\Label\ConnectionForm$$
Then $\Tt{\lambda}{f}:=F_{\lambda}e_1$ is called $T_{\!\lambda}$-transform, or spectral
transform of $f$.} \TPic

\noindent
This way, the $T_{\!\lambda}$-transform $\Tt{\lambda}{f}$ of an isothermic immersion
$f$ is unique up to M\"obius transformation.
On the other hand, if $\tilde f$ is an essential M\"obius transform of $f$,
say $\tilde f=-(f-m)^{-1}$, then $\Ct{\tilde{f}}=G\Ct{f}$ is a Goursat transform
of $\Ct{f}$ and $d\tilde{F}_{\lambda}=\tilde{F}_{\lambda}\tilde{\Phi}_{\lambda}$
with $\tilde{F}_{\lambda}=F_{\lambda}\cdot\left(\matrix{(f-m)^{-1}&-1\cr0&(f-m)
\cr}\right).$
Consequently, $\tilde{F}_{\lambda}e_1$ and $F_{\lambda}e_1$ are just
different homogeneous coordinates for the same surface in $S^3$: the
$T_{\!\lambda}$-transformation is a well defined transformation for
isothermic surfaces in M\"obius geometry.
In fact, the $T_{\!\lambda}$-transform of an isothermic immersion
can be defined using {\it any} adapted frame: observing the effect of
a Gauge transformation of adapted frames, $F_0\to\tilde{F}_0=F_0A$ with
$Ae_1=e_1a$ and $a:M^2\to\H$, on the corresponding connection forms shows
that a suitable definition for $\tilde{\Phi}_{\lambda}$ provides a
family of integrable connection forms corresponding to \ConnectionForm.
Given an adapted frame $\tilde{F}_0:M^2\to Sl(2,\H)$ of an isothermic
immersion $f:(M^2,q)\to{\rm Im}\H$,
$$\tilde{F}_0^{-1}d\tilde{F}_0=\tilde{\Phi}_0
	= \left(\matrix{\varphi&\hat{\psi}\cr\psi&\hat{\varphi}\cr}\right),
	\hskip1em
	\tilde{\Phi}_{\lambda}
	:= \left(\matrix{\varphi&\hat{\psi}+\lambda\psi^{\ast}
	\cr\psi&\hat{\varphi}\cr}\right)\Eqno\Label\Independent$$
is integrable with $\tilde{F}_{\lambda}e_1\simeq\Tt{\lambda}{f}$.
Herein, $\psi^{\ast}$ is defined --- analogously to $d\Ct{f}$ --- by the
equation\Footnote{Note that $\psi^{\perp}$ (similar to $df^{\perp}$) defines
a complex line bundle over $M^2$ as $(\tilde{F}_0e_1)\H:M^2\to S^3$.}
$\psi\cdot\psi^{\ast}\simeq q$.

To understand the geometry of $\Tt{\lambda}{f}$, we first note that
$\Phi_{\lambda}$ defined by \ConnectionForm\ takes values in the Lie
algebra of the M\"obius group of $S^3\subset\HP^1$. Consequently,
$\Tt{\lambda}{f}$ can be assumed (after a suitable M\"obius transformation)
to take values in $S^3$, $\Tt{\lambda}{f}:M^2\to S^3={\rm Im}\H\cup\{\infty\}$.
Using affine coordinates, $F_{\lambda}=\left(\matrix{gb&\hat{g}\hat{b}\cr
b&\hat{b}\cr}\right)$, yields $dg=(\hat{g}-g)\hat{b}\,df\,b^{-1}$ and
$d\hat{g}=(g-\hat{g})b\,\lambda d\Ct{f}\,\hat{b}^{-1}$, hence $dg\wedge
(g-\hat{g})^{-1}d\hat{g}=0$.
According to the above lemma, $g\simeq\Tt{\lambda}{f}$ is isothermic with
Christoffel transform $\Ct{g}$ given by $\mu d\Ct{g}=(\hat{g}-g)^{-1}d\hat{g}
(\hat{g}-g)^{-1}$.
Calculating the corresponding quadratic differential, $\mu dg\cdot d\Ct{g}
=-\lambda[(g-\hat{g})\hat{b}]\,df\cdot d\Ct{f}[(g-\hat{g})\hat{b}]^{-1}
\simeq-\lambda q$, shows that the curvature lines of $f$ and $\Tt{\lambda}{f}$
correspond\Footnote{In fact, the family $\lambda\to\Tt{\lambda}{f}$ is the
``conformal deformation'' \MussoDef\ of $f$ (cf.\Cartan): isothermic surfaces
are the only non rigid surfaces in M\"obius geometry with prescribed
curvature lines and conformal metric (Calapso potential), resp.\ with
prescribed ``conformal Hopf differential'' \Joerg.}.
Moreover, the Christoffel transform $\Ct{g}$ of $g$ can be scaled,
$\mu=-\lambda$, such that $dg\cdot d\Ct{g}\simeq q$.
Then, $\hat{g}=\Dt{-\lambda}{g}$;
on the other hand, $\hat{g}=\Tt{\lambda}{\Ct{f}}$ as the Gauge transformation
$(e_1,e_2)\to({1\over\lambda}\sqrt{|\lambda|}e_2,\sqrt{|\lambda|}e_1)$
interchanges the roles of $df$ and $d\Ct{f}$ in \ConnectionForm, i.e.\
$\hat{g}\simeq\Tt{\lambda}{\Ct{f}}$.
This way, we obtain the first

\SubSection{Permutability}
theorem (cf.\Bianchi), $\Tt{\lambda}{\Ct{f}}=\Dt{-\lambda}{\Tt{\lambda}{f}}$.
As the $\Tt{\lambda}{}$-transformation can be defined using {\it any} adapted
frame, it is obvious that the 1-parameter family $\lambda\to\Tt{\lambda}{f}$
of $T$-transforms of an isothermic immersion $f$ is a 1-parameter {\it group}:
$\Tt{\lambda}{\Tt{\mu}{f}}=\Tt{\lambda+\mu}{f}$.
Therefore, $\Tt{\mu}{\Dt{\lambda}{f}}=\Dt{\lambda-\mu}{\Tt{\mu}{f}}$.

Given an isothermic immersion $f:(M^2,q)\to{\rm Im}\H$, and a Darboux
transform $\Dt{\lambda}{f}$ of $f$, $$\matrix{
   \lambda d\Ct{f} &=&
    (\Dt{\lambda}{f}-f)^{-1} d\Dt{\lambda}{f} (\Dt{\lambda}{f}-f)^{-1},\cr
   \lambda d\Ct{\Dt{\lambda}{f}} &=&
    (\Dt{\lambda}{f}-f)^{-1} df (\Dt{\lambda}{f}-f)^{-1}.\cr}$$
Hence, $\Ct{f}$ and $\Ct{\Dt{\lambda}{f}}$ can be positioned such that
$\lambda(\Ct{\Dt{\lambda}{f}}-\Ct{f})=(\Dt{\lambda}{f}-f)^{-1}$.
Then, $d(\Ct{\Dt{\lambda}{f}}-\Ct{f})=(\Ct{\Dt{\lambda}{f}}-\Ct{f})\lambda
df(\Ct{\Dt{\lambda}{f}}-\Ct{f})-d\Ct{f}$ showing (cf.\Bianchi,\JePe) that
$\Ct{\Dt{\lambda}{f}}=\Dt{\lambda}{\Ct{f}}$.
Combining this permutability theorem with the first, we obtain the following
scheme:
\Figure{\epsfbox{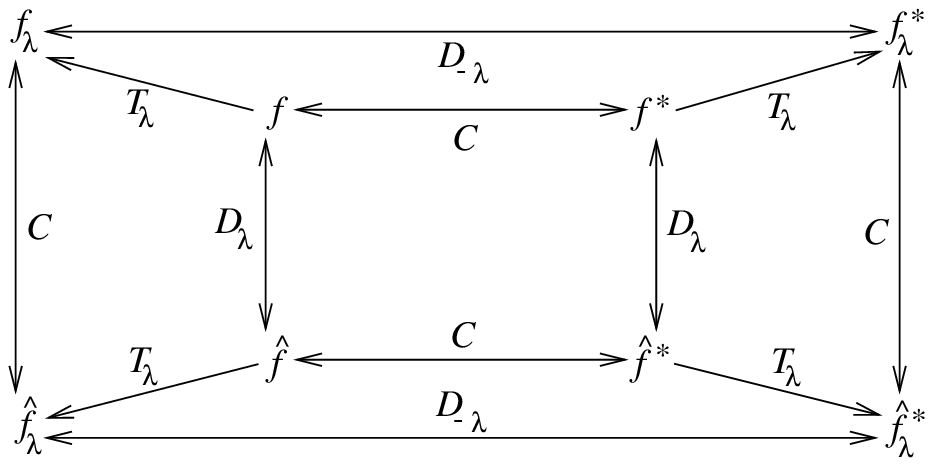}}{A permutability theorem}\Label\PermThm

\noindent
Finally, consider two $\Dt{\lambda}{}$-transforms $(\Dt{\lambda}{f})_1$ and
$(\Dt{\lambda}{f})_2$ of an isothermic immersion $f:(M^2,q)\to{\rm Im}\H$.
The $\Tt{\lambda}{}$-transforms of these three surfaces are only well defined
up to M\"obius transformation. But, after suitable M\"obius transformations,
$\Tt{\lambda}{(\Dt{\lambda}{f})_1}$ is the Christoffel transform of
$\Tt{\lambda}{f}$, as well as $\Tt{\lambda}{(\Dt{\lambda}{f})_2}$ is
--- after possibly other M\"obius transformations. Thus,
$\Tt{\lambda}{(\Dt{\lambda}{f})_2}=\Gt{\Tt{\lambda}{(\Dt{\lambda}{f})_1}}$:
after a $\Tt{\lambda}{}$-transformation, the ``difference'' between two
$\Dt{\lambda}{}$-transforms of an isothermic surface is given by a Goursat
transformation.

\Section{Surfaces of constant mean curvature 1 in $H^3$}
Thus far, we have not yet given any examples of Darboux or $T$-transforms.
To do so, it turns out to be convenient to use a suitably adapted frame for
the two envelopes $f,\hat{f}:M^2\to S^3\simeq s_1$ of a sphere congruence
$s:M^2\to\R^6_1$: given that the curvature lines on both envelopes do
correspond ($s$ is ``Ribaucour'', $\hat{f}$ is a ``Ribaucour transformation''
of $f$), there are two congruences of 2-spheres in $S^3$ such that each one
intersects both surfaces orthogonally to one family of curvature lines.
Moreover, homogeneous coordinates $f,\hat{f}:M^2\to\H^2$ can be chosen
such that these sphere congruences become  $s_j$ and $s_k$ with respect to
the frame $F=(f,\hat{f}):M^2\to{\rm Sl}(2,\H)$, $s=s_i$ (cf.\ClassicalFrameS),
and that $\langle df,\hat{f}\rangle=0$.
This way, the homogeneous coordinates of the two envelopes are fixed
up to constant rescaling. We will refer to such a frame as a ``Ribaucour
frame'' for $f$ and $\hat{f}$.

If we additionally assume $f:(M^2,q)\to S^3$ to be isothermic, we can introduce
conformal curvature line parameters $z:M^2\to\C$. Then, the connection form of
a Ribaucour frame, $\Phi=F^{-1}dF$, takes the form $$\Phi=\left(\matrix{
	{i\over2}[\ast du-(He^udz+\hat{H}e^{-u}d\bar{z})j]&
	(\hat{\lambda}e^udz-\lambda e^{-u}d\bar{z})j\cr
	e^udz\,j&
	{i\over2}[\ast du-(He^udz+\hat{H}e^{-u}d\bar{z})j]\cr
	}\right)\Eqno\Label\RibConnectionForm$$
with suitable real functions $u,H,\hat{H},\lambda,\hat{\lambda}:M^2\to\R$.
The Maurer-Cartan equation for $\Phi$, $0=d\Phi+\Phi\wedge\Phi$, yields
the Gauss equation $$\Delta u = -H^2e^{2u}+\hat{H}^2e^{-2u} \,+
	4e^{2u}\hat{\lambda},\Eqno\Label\Gauss$$ and the Codazzi equations
$$\matrix{0=\hat{\lambda}_{\bar{z}}e^u+\lambda_ze^{-u} && {\rm and} &&
	H_{\bar{z}}e^u=\hat{H}_ze^{-u} \cr }.\Eqno\Label\Codazzi$$

\SubSection{The D-transformation and cmc-1 surfaces in $H^3$}
Now, $f$ and $\hat{f}$ induce conformally equivalent metrics on $M^2$
if either $\lambda\equiv0$ or $\hat{\lambda}\equiv0$ while the other is
constant by \Codazzi. In the former case, the two surfaces, $f$ and $\hat{f}$,
are M\"obius equivalent; in the latter $\hat{f}$ is a Darboux transform of
$f$ (and vice versa), and $F$ is a ``curved flat frame'' (cf.\Suppl).

The second surface, $\hat{f}$, of such a pair $F=(f,\hat{f})$ of Darboux
transforms of each other is totally umbilic if and only if $H\equiv0$ in
\RibConnectionForm.
Obviously, this is equivalent to the fact that ${\rm tr}\langle df,ds\rangle
=0$ with respect to the induced conformal structure $\langle df,df\rangle$
of $f$, characterizing the enveloped sphere congruence $s$ as the ``central
sphere congruence'' of $f$ (cf.\Blaschke).
By the Codazzi equations, also $\hat{H}$ is constant; consequently,
$\hat{H}\hat{f}+\lambda\,s=:s_{\infty}$ is constant, $ds_{\infty}=
\hat{H}d\hat{f}+\lambda\,ds\equiv0$, thus describing a constant 2-sphere
(as long as $\lambda\neq0$).
As $\langle\hat{f},s_{\infty}\rangle\equiv0$, $\hat{f}$ takes values in
that 2-sphere $s_{\infty}$. 
If $\hat{H}\equiv0$, also $f$ takes values in $s_{\infty}$ and a ``Darboux
pair of meromorphic functions'' is obtained as a degenerate case.
If $\hat{H}\neq0$, we can assume $\hat{H}\equiv1$ without loss of
generality\Footnote{By a constant Gauge transformation $(f,\hat{f})\to
({1\over\hat{H}}f,\hat{H}\hat{f})$ and possibly interchanging the roles of the
principal curvature directions, a new Ribaucour frame with \RibConnectionForm\
is obtained, where $u-{\rm ln}|\hat{H}|\to u$.}\Label\Loss.
Interpreting the two components of $$\matrix{
S^3\setminus s_{\infty}\cong\{x\in\R^6_1\,|\,\langle x,x\rangle=0,\,
	\langle x,s_1\rangle=0,\,\langle x,s_{\infty}\rangle=-{1\over2}\}}$$
as hyperbolic spaces of sectional curvature $-4\lambda^2$ with infinity
boundary $s_{\infty}$, the sphere congruence $s$ becomes a congruence
of horospheres in one of these hyperbolic spaces, $H^3\subset S^3$.
Using the tangent plane congruence $t:=s+2\lambda f$ of $f$ in hyperbolic
space\Footnote{As $\langle t,s_{\infty}\rangle\equiv0$, the spheres $t$
intersect the infinity boundary $s_{\infty}$ orthogonally.}, the first
and second fundamental form of $f$ become $$\matrix{
	\hfill\langle df,df\rangle  &=&\hfill e^{2u}\,(dx^2+dy^2)&&\cr
	-\langle df,dt\rangle &=& -2\lambda e^{2u}\,(dx^2+dy^2)&+&(dx^2-dy^2)\cr
	}\Eqno\Label\FundamentalForms$$
showing that $f$ has constant mean curvature\Footnote{Here, we could as
well have argued that the mean curvature of $f$ at any point coincides with
the mean curvature $-2\langle s,s_{\infty}\rangle$ of the central sphere
$s$ at that point, measured in any space form.} $h\equiv-2\lambda$ (as a
surface in hyperbolic space).
Summarizing, we obtain the following

\Theorem{}{Let $\hat{f}:(M^2,\bar{q})\to s_{\infty}\subset S^3$ be an
isothermic immersion into a 2-sphere $s_{\infty}$. Then, any Darboux transform
$f:(M^2,q)\to S^3$ of $\hat{f}$ either takes values in $s_{\infty}$, too,
or parametrizes a surface of constant mean curvature $h$ in a hyperbolic
space $H^3_k$ of sectional curvature $k=-h^2$.
In the latter case, the enveloped Ribaucour sphere congruence is the central
sphere congruence of $f$, and $\hat{f}$ is the hyperbolic Gauss map of $f$.}

\noindent
To confirm the last statement of the theorem, we just notice that the circles
intersecting the central spheres $s$ orthogonally in $f$ and $\hat{f}$ are
the hyperbolic geodesics orthogonal to $f$.
Their intersection points with the infinity boundary $s_{\infty}$ are given
by $\hat{f}$ and $\tilde{f}:=f+{1\over2\lambda}s+{1\over4\lambda^2}\hat{f}$.
As $f$ and $\hat{f}$ induce conformally equivalent metrics on $M^2$ (while
$f$ and $\tilde{f}$ do not) the claim follows (see \Bryant).

\SubSection{The T-transformation and cmc-1 surfaces in $H^3$}
Having a closer look at the connection form $$\Phi_{\lambda}:=\left(\matrix{
	{i\over2}[\ast du-e^{-u}d\bar{z}\,j]&-\lambda e^{-u}d\bar{z}\,j\cr
	e^udz\,j&{i\over2}[\ast du-e^{-u}d\bar{z}\,j]\cr
	}\right)\Eqno\Label\CMCConnectionForm$$
of the above Ribaucour frame for the surface $f_{\lambda}:=f:M^2\to
H^3_{\-4\lambda^2}$ of constant mean curvature $-2\lambda$ and its hyperbolic
Gauss map $\hat{f}_{\lambda}:=\hat{f}:M^2\to s_{\infty}=:s_{\infty}(\lambda)$
shows immediately\Footnote{Remember that the $T$-transformation is frame
independent, see \Independent.} that $f_{\lambda}=\Tt{\lambda}{f_0}$ is a
$\Tt{\lambda}{}$-transform of an isothermic surface $f_0:(M^2,q)\to{\rm Im}\H$:
sending $\lambda\to0$, the infinity sphere $s_{\infty}(\lambda)\to\hat{f}_0
\equiv p_{\infty}$ degenerates to a point, $d\hat{f}_0\equiv0$, that we
interprete as the point at infinity of Euclidean 3-space ${\rm Im}\H$.
Since the central spheres $s_0$ of $f_0$ all contain the point $p_{\infty}$
at infinity they are planes in ${\rm Im}\H$. Consequently, the surface
$f_0$ is minimal in ${\rm Im}\H$.
On the other hand, away from umbilics, every minimal immersion $f_0:M^2\to{\rm
Im}\H$ carries conformal curvature line coordinates $z:M^2\to\C$ (minimal
immersions are isothermic) such that the Hopf differential $q=dz^2$.
Thus, every minimal immersion has a frame with connection form $\Phi_0$, and

\Theorem{}{The surfaces of constant mean curvature $h=-2\lambda$
in hyperbolic space of sectional curvature $k=-4\lambda^2$ are in
one-to-one correspondence with the minimal surfaces in Euclidean
3-space via the $\Tt{\lambda}{}$-transformation.}

\noindent
As all surfaces $f_{\lambda}=\Tt{\lambda}{f_0}:M^2\to H^3_{-4\lambda^2}$,
$\lambda\in\R$, are isometric, $I_{\lambda}=I_0$, and their second
fundamental forms satisfy $\II_{\lambda}=\II_0-2\lambda I_0$ the
$\Tt{\lambda}{}$-transformation yields exactly the ``Umehara-Yamada
perturbation'' \UmYa\ of minimal surfaces into constant mean curvature
surfaces in hyperbolic space.
Consequently, our previously discussed version \WeierstrassRepresentation\
of the Weierstrass representation for minimal surfaces in Euclidean 3-space
``perturbs'' into a version of Bryant's Weierstrass type representation \Bryant:
given Weierstrass data $(g,\omega)$ and $\lambda\in\R$ the corresponding surface
$f_{\lambda}$ of constant mean curvature $-2\lambda$ in $H^3_{-4\lambda^2}$
(together with its hyperbolic Gauss map $\hat{f}_{\lambda}$, if $\lambda\neq0$)
or $\R^3$, respectively, is obtained by integrating the system $$\matrix{
	df_{\lambda} &=& \hat{f}_{\lambda}\cdot
		[{1\over2}(i-jg)\omega\,j(i-jg)],\hfill\cr
	d\hat{f}_{\lambda} &=& f_{\lambda}\cdot
		[-2\lambda(i-jg)^{-1}jdg(i-jg)^{-1}].\cr}
	\Eqno\Label\BryantRepresentation$$

\Remark{The Gauss map $(i-jg)i(i-jg)^{-1}$ of the minimal surface $f_0$ is called the
``secondary Gauss map'' of $f_{\lambda}$, and $f_0$ its ``minimal cousin'' (cf.\Bryant).}

\SubSection{Isothermic surfaces of spherical type}
The Gauss equation \Gauss\ for the connection form \CMCConnectionForm\
associated with a surface $f_{\lambda}:M^2\to H^3_{-4\lambda^2}$ of
constant mean curvature $-2\lambda$ (or, its minimal cousin $f_0:M^2\to\R^3$)
reduces to the Liouville equation
$$\Delta u = e^{-2u}.\Eqno\Label\Liouville$$
Consequently, the induced metric $\langle ds,ds\rangle=e^{-2u}(dx^2+dy^2)$
of the (central) sphere congruence $s$ has constant Gauss curvature $1$:

\Definition{~(spherical type)}{An isothermic immersion is called isothermic
of spherical type if its central sphere congruence induces a metric of constant
Gauss curvature 1.}

\noindent
We will see that this is a (M\"obius geometric) characterization for surfaces
of constant mean curvature $h$ in hyperbolic space of sectional curvature
$k=-h^2$, and minimal surfaces in Euclidean space:
%
it is well known \Blaschke\ that a surface $f:M^2\to S^3$ is isothermic if
and only if its central sphere congruence $s$ is Ribaucour.
Therefore, we can associate a Ribaucour frame $F=(f,\hat{f})$ to $f$ where
$\hat{f}$ is the second envelope of the central sphere congruence.
In the connection form \RibConnectionForm, the fact that $s=s_i$ is the
central sphere congruence of $f$ is reflected by $H\equiv0$.
By \Codazzi, $\hat{H}$ is constant so that without loss of generality
$\hat{H}\equiv1$ (see footnote \Loss). 
Now, $u$ satisfies the Liouville equation \Liouville\
if and only if $\hat{\lambda}\equiv0$,
if and only if the connection form \RibConnectionForm\ takes the form
\CMCConnectionForm. Thus,

\Theorem{}{An isothermic immersion is isothermic of spherical type\newline
iff the second envelope of its central sphere congruence is a Darboux
transform \centerline{--- this Darboux transform is then either totally
umbilic, or a point ---}\newline
iff the surface is either a surface of constant mean curvature
$h$ in a hyperbolic space of sectional curvature $k=-h^2$, or a minimal
surface in Euclidean space.}

\noindent
Consequently, {\it every} surface $f:M^2\to H^3_k$ of constant mean
curvature $h$ in a space of sectional curvature $k=-h^2$ can be obtained
as a $\Dt{\lambda}{}$-transform, $\lambda=-{h\over2}$, of a totally umbilic
immersion $n_h:M^2\to S^2\cong\partial H^3_k$.
This gives rise to another Weierstrass type representation for such surfaces:
given Weierstrass data $(g,\omega)$, integration of Darboux's linear system%
\Footnote{Note that multiplying the second equation by $j$ and substituting
$jv_2\to v_2$, this system can entirely be handled in the context of complex
function theory: given a (complex) fundamental system {\it any} solution is
obtained as a {\it quaternionic} superposition.}
$$\matrix{dv_1 &+& \lambda\omega j\cdot v_2 &=& 0 \cr
	dv_2 &-& jdg\cdot v_1 &=& 0\cr}\Eqno\Label\Darboux$$
yields a surface $f_{\lambda}^{\#}=-jg+v_2v_1^{-1}$ of constant mean curvature
$h=2\lambda$ in $H^3_{-4\lambda^2}$ as a $\Dt{\lambda}{}$-transform of its
hyperbolic Gauss map $n_h=-jg$ --- as long as an initial condition for $v$
is chosen such that $v_2v_1^{-1}\not\in\C j$.
The secondary Gauss map of $f_{\lambda}^{\#}$ and its minimal cousin are
obtained as $\Tt{\lambda}{}$-transforms of $-jg$ and $f_{\lambda}^{\#}$,
respectively.
Thus, comparing this second Weierstrass type representation for surfaces of
constant mean curvature $h$ in $H^3_{-h^2}$ with Bryant's, the roles of the
hyperbolic and secondary Gauss maps of the surface are interchanged.

\SubSection{Duality for cmc-1 surfaces in $H^3$}
In \UmYaDual\ M.~Umehara and K.~Yamada introduced a notion of ``duality'' for
cmc-1 surfaces in hyperbolic space: given a cmc-1 surface $f:M^2\to H^3$ with
hyperbolic and secondary Gauss maps $n_h,n_s:M^2\to S^2\cong\partial H^3$, the
``dual cmc-1 surface'' $f^{\#}:M^2\to H^3$ is obtained by interchanging the
roles of the two Gauss maps, $n_h^{\#}=n_s$ and $n_s^{\#}=n_h$.
However, the hyperbolic Gauss map $n_h$ is only well defined up to M\"obius
transformation of $S^2$ since isometries of $H^3$ extend to M\"obius
transformations of $S^2\cong\partial H^3$.
Consequently, the dual cmc-1 surface $f^{\#}$ depends on the position of $f$
in hyperbolic space; as there is a 3-parameter family of essential M\"obius
transformations of the 2-sphere, generally, there is a 3-parameter family of
dual cmc-1 surfaces.
In the context of transformations of isothermic surfaces, this duality
relation occurs as a special case of the permutability theorem sketched
in Fig.~\PermThm: as discussed above, the hyperbolic Gauss map $n_h=\hat{f}$ 
of a cmc-1 surface $f$ is a $\Dt{1\over2}{}$-transform of $f$ where the
enveloped sphere congruence is the central sphere congruence of $f$.
The minimal cousin of $f$ and its Gauss map (the secondary Gauss map
of $f$) are obtained as $\Tt{1\over2}{}$-transforms of $f$ and $\hat{f}$,
respectively.
Now, the dual cmc-1 surface\Footnote{With respect to the original polarization,
$f^{\#}$ has constant mean curvature -1; changing the direction of the
normal field, $q\to-q$, yields mean curvature 1 (cf.\UmYaDual).}
$f^{\#}=\hat{f}^{\ast}_{1\over2}$ is a $\Dt{-{1\over2}}{}$-transform of
$n_s=\Tt{1\over2}{n_h}$, and $f^{\ast}=\Tt{1\over2}{f^{\#}}$ is its minimal
cousin.
The 3-parameter ambiguity of the dual cmc-1 surface $f^{\#}$ is the
3-parameter ambiguity of the Darboux transformation, and the minimal cousins
of different dual cmc-1 surfaces differ by a Goursat transformation.

We conclude our paper with a simple
\SubSection{Example:}
consider $f(z):=-jz$ as an isothermic immersion of the polarized plane $(\C,d\bar{z}^2)$
into $\C j\subset{\rm Im}\H$.
Its Christoffel transform is $\Ct{f}=zj:(\C,dz^2)\to\C j$.
With $$ \Phi_{\lambda}
	= \left(\matrix{0&\lambda d\Ct{f}\cr df&0\cr}\right)
	= \left(\matrix{1&0\cr0&-j\cr}\right)
	\left(\matrix{0&\lambda dz\cr dz&0\cr}\right)
	\left(\matrix{1&0\cr0&j\cr}\right), $$
integration of the linear system $dF_{\lambda}=F_{\lambda}\Phi_{\lambda}$,
$F_{\lambda}(0)=Id$, yields $$ F_{\lambda} =
	\left(\matrix{1&0\cr0&-j\cr}\right)
	\left(\matrix{\cosh(\sqrt{\lambda}z)&
		\sqrt{\lambda}\sinh(\sqrt{\lambda}z)\cr
		{1\over\sqrt{\lambda}}\sinh(\sqrt{\lambda}z)&
		\cosh(\sqrt{\lambda}z)\cr}\right)
	\left(\matrix{1&0\cr0&j\cr}\right). $$
Thus, $\Tt{\lambda}{f}=-j{1\over\sqrt{\lambda}}\tanh(\sqrt{\lambda}z)$ --- or
any M\"obius transform of it. Its Christoffel transform $\Ct{\Tt{\lambda}{f}}$
depends on the choice of representative for $\Tt{\lambda}{f}$: with the above
choice, $\Ct{\Tt{\lambda}{f}}={1\over2}[z+{1\over2\sqrt{\lambda}}\sinh(2\sqrt{
\lambda}z)]j$, any other choice yields a Goursat transform thereof.
For example, stereographic projection of $\Tt{\lambda}{f}$ to the unit
2-sphere yields a family
${1\over4}\{{\rm Re}[{\cosh(2\sqrt{\lambda}z)-1\over\lambda}]\,i
+[z+{\sinh(2\sqrt{\lambda}z)\over2\sqrt{\lambda}}]\,j
+j\,{1\over\lambda}[z-{\sinh(2\sqrt{\lambda}z)\over2\sqrt{\lambda}}]\}$
of minimal surfaces with the Catenoid at $\lambda=1$ and the Enneper
surface at $\lambda=0$.
The homogeneous coordinates of any Darboux transform $\Dt{\lambda}{f}=
F_0F_{\lambda}^{-1}v_0$ with constant $v_0\in\H^2$.
Choosing $\Dt{\lambda}{f}(0)\simeq v_0\H\not\in\C j\cup\{\infty\}$,
$\Dt{\lambda}{f}$ is a surface of constant mean curvature $2\lambda$
in $H^3_{-4\lambda^2}$.
For example, $v_0=(1,-i)^t$ yields the surfaces $\Dt{\lambda}{f}(z)=-j\{z-
[{\sinh(\sqrt{\lambda}z)\over\sqrt{\lambda}}-\cosh(\sqrt{\lambda}z)k]
[\cosh(\sqrt{\lambda}z)-\sqrt{\lambda}\sinh(\sqrt{\lambda}z)k]^{-1}\}$
that have the above minimal surfaces as their minimal cousins.

In order to determine the Darboux transforms of $\Tt{\lambda}{f}$ too,
one has to integrate Darboux's linear system again with a connection form
$\Phi_{\lambda,\mu}=\Phi_{\mu}(\Tt{\lambda}{f})$.
As the Darboux transformation is a M\"obius geometric notion any choice
of representative for $\Tt{\lambda}{f}$ will lead to the same (in terms
of M\"obius geometry) result.  Thus, with the above choice
$\Tt{\lambda}{f}=-j{\tanh(\sqrt{\lambda}z)\over\sqrt{\lambda}}$, we obtain
$F_{\lambda,\mu}e_1=F_{\lambda+\mu}{e_1\over\cosh(\sqrt{\lambda}z)}$
confirming that $\Tt{\mu}{\Tt{\lambda}{f}}=\Tt{\lambda+\mu}{f}$.
In the special case $\mu=-\lambda$ of our permutability theorem (see Figure
\PermThm), $F_{\lambda,0}F_{\lambda,-\lambda}^{-1}=[F_0F_{\lambda}^{-1}]^{-1}$
so that the Darboux transforms $\Dt{-\lambda}{\Tt{\lambda}{f}}=-j\{
{\tanh(\sqrt{\lambda}z)\over\sqrt{\lambda}}-{1\over\cosh(\sqrt{\lambda}z)}
[z-k][\cosh(\sqrt{\lambda}z)-\sqrt{\lambda}\sinh(\sqrt{\lambda}z)(z-k)]^{-1}\}$
can be determined without further integration.
Note that according to the permutability theorem all the surfaces
$\Dt{-\lambda}{\Tt{\lambda}{f}}$ have the Enneper surface as their
minimal cousin,
and they are dual cmc surfaces of $\Dt{\lambda}{f}$ in the sense of
\UmYaDual.

\Acknowledgements{We would like to thank Fran Burstall for many fruitful
discussions and comments.

Also, we would like to thank the Department of Mathematics and the GANG
at the University of Massachusetts in Amherst as well as the Department
of Mathematics at the University of l'Aquila for their hospitality when
the third and the first author, respectively, were visiting.}

\References

\bye